\documentclass[english,envcountsame,envcountsect]{svmult}

\usepackage{babel}
\usepackage{amstext}
\usepackage{amsmath}
\usepackage{amsfonts}
\usepackage{ifthen}
\usepackage{xypic}
\usepackage{url}
\usepackage{makeidx}         
\usepackage{graphicx}        
\usepackage{multicol}        
\usepackage[bottom]{footmisc}
\newcommand{\myqed}{\hspace*{\fill}\qed}
\spnewtheorem{properties}[theorem]{Properties}{\bfseries}{\rmfamily}
\spnewtheorem{assumption}[theorem]{Assumption}{\bfseries}{\rmfamily}
\catcode`\@=11
\def\opn#1#2{\def#1{\mathop{\kern0pt\fam0#2}\nolimits}}
\def\underrightarrow{\mathpalette\underrightarrow@}
\def\underrightarrow@#1#2{\vtop{\ialign{$##$\cr
 \hfil#1#2\hfil\cr\noalign{\nointerlineskip}%
 #1{-}\mkern-6mu\cleaders\hbox{$#1\mkern-2mu{-}\mkern-2mu$}\hfill
 \mkern-6mu{\to}\cr}}}

\def\underleftarrow{\mathpalette\underleftarrow@}
\def\underleftarrow@#1#2{\vtop{\ialign{$##$\cr
 \hfil#1#2\hfil\cr\noalign{\nointerlineskip}#1{\leftarrow}\mkern-6mu
 \cleaders\hbox{$#1\mkern-2mu{-}\mkern-2mu$}\hfill
 \mkern-6mu{-}\cr}}}
\let\amp@rs@nd@\relax
\newdimen\ex@
\newdimen\bigaw@
\newdimen\minaw@
\newdimen\minCDaw@
\newif\ifCD@
\def\minCDarrowwidth#1{\minCDaw@#1}
\newenvironment{CD}{\@CD}{\@endCD}
\def\@CD{\def\A##1A##2A{\llap{$\vcenter{\hbox
 {$\scriptstyle##1$}}$}\Big\uparrow\rlap{$\vcenter{\hbox{%
$\scriptstyle##2$}}$}&&}%
\def\V##1V##2V{\llap{$\vcenter{\hbox
 {$\scriptstyle##1$}}$}\Big\downarrow\rlap{$\vcenter{\hbox{%
$\scriptstyle##2$}}$}&&}%
\def\={&\hskip.5em\mathrel
 {\vbox{\hrule width\minCDaw@\vskip3\ex@\hrule width
 \minCDaw@}}\hskip.5em&}%
\def\verteq{\Big\Vert&&}%
\def\noarr{&&}%
\def\vspace##1{\noalign{\vskip##1\relax}}\relax\let\amp@rs@nd@&\iffalse}\fi
 \CD@true\vcenter\bgroup\relax\let\\=\cr\iffalse}\fi\tabskip\z@skip\baselineskip20\ex@
 &\hfill$\m@th##$\hfill\cr}
\def\@endCD{\cr\egroup\egroup}
\def\>#1>#2>{\amp@rs@nd@\setbox\z@\hbox{$\scriptstyle
 \;{#1}\;\;$}\setbox\@ne\hbox{$\scriptstyle\;{#2}\;\;$}\setbox\tw@
 \hbox{$#2$}\ifCD@
 \global\bigaw@\minCDaw@\else\global\bigaw@\minaw@\fi
 \ifdim\wd\z@>\bigaw@\global\bigaw@\wd\z@\fi
 \ifdim\wd\@ne>\bigaw@\global\bigaw@\wd\@ne\fi
 \ifCD@\hskip.5em\fi
 \ifdim\wd\tw@>\z@
 \mathrel{\mathop{\hbox to\bigaw@{\rightarrowfill}}\limits^{#1}_{#2}}\else
 \mathrel{\mathop{\hbox to\bigaw@{\rightarrowfill}}\limits^{#1}}\fi
 \ifCD@\hskip.5em\fi\amp@rs@nd@}
\def\<#1<#2<{\amp@rs@nd@\setbox\z@\hbox{$\scriptstyle
 \;\;{#1}\;$}\setbox\@ne\hbox{$\scriptstyle\;\;{#2}\;$}\setbox\tw@
 \hbox{$#2$}\ifCD@
 \global\bigaw@\minCDaw@\else\global\bigaw@\minaw@\fi
 \ifdim\wd\z@>\bigaw@\global\bigaw@\wd\z@\fi
 \ifdim\wd\@ne>\bigaw@\global\bigaw@\wd\@ne\fi
 \ifCD@\hskip.5em\fi
 \ifdim\wd\tw@>\z@
 \mathrel{\mathop{\hbox to\bigaw@{\leftarrowfill}}\limits^{#1}_{#2}}\else
 \mathrel{\mathop{\hbox to\bigaw@{\leftarrowfill}}\limits^{#1}}\fi
 \ifCD@\hskip.5em\fi\amp@rs@nd@}

\def\@CDS{\def\A##1A##2A{\llap{$\vcenter{\hbox
 {$\scriptstyle##1$}}$}\Big\uparrow\rlap{$\vcenter{\hbox{%
$\scriptstyle##2$}}$}&}%
\def\V##1V##2V{\llap{$\vcenter{\hbox
 {$\scriptstyle##1$}}$}\Big\downarrow\rlap{$\vcenter{\hbox{%
$\scriptstyle##2$}}$}&}%
\def\={&\hskip.5em\mathrel
 {\vbox{\hrule width\minCDaw@\vskip3\ex@\hrule width
 \minCDaw@}}\hskip.5em&}
\def\verteq{\Big\Vert&}
\def\novarr{&}
\def\noharr{&&}
\def\SE##1E##2E{\slantedarrow(0,18)(4,-3){##1}{##2}&}
\def\SW##1W##2W{\slantedarrow(24,18)(-4,-3){##1}{##2}&}
\def\NE##1E##2E{\slantedarrow(0,0)(4,3){##1}{##2}&}
\def\NW##1W##2W{\slantedarrow(24,0)(-4,3){##1}{##2}&}
\def\slantedarrow(##1)(##2)##3##4{%
\thinlines\unitlength1pt\lower 6.5pt\hbox{\begin{picture}(24,18)%
\put(##1){\vector(##2){24}}%
\put(0,8){$\scriptstyle##3$}%
\put(20,8){$\scriptstyle##4$}%
\end{picture}}}
\def\vspace##1{\noalign{\vskip##1\relax}}\relax\let\amp@rs@nd@&\iffalse}\fi
 \CD@true\vcenter\bgroup\relax\let\\=\cr\iffalse}\fi\tabskip\z@skip\baselineskip20\ex@
 &\hfill$\m@th##$\hfill\cr}
\def\@endCDS{\cr\egroup\egroup}
\newdimen\TriCDarrw@
\newif\ifTriV@

\def\@TriCDV{\TriV@true\def\TriCDpos@{6}\@TriCD}
\def\@TriCDA{\TriV@false\def\TriCDpos@{10}\@TriCD}
\def\@TriCD#1#2#3#4#5#6{%
\setbox0\hbox{$\ifTriV@#6\else#1\fi$}
\TriCDarrw@=\wd0 \advance\TriCDarrw@ 24pt
\advance\TriCDarrw@ -1em
\def\SE##1E##2E{\slantedarrow(0,18)(2,-3){##1}{##2}&}
\def\SW##1W##2W{\slantedarrow(12,18)(-2,-3){##1}{##2}&}
\def\NE##1E##2E{\slantedarrow(0,0)(2,3){##1}{##2}&}
\def\NW##1W##2W{\slantedarrow(12,0)(-2,3){##1}{##2}&}
\def\slantedarrow(##1)(##2)##3##4{\thinlines\unitlength1pt
\lower 6.5pt\hbox{\begin{picture}(12,18)%
\put(##1){\vector(##2){12}}%
\put(-4,\TriCDpos@){$\scriptstyle##3$}%
\put(12,\TriCDpos@){$\scriptstyle##4$}%
\end{picture}}}
\def\={\mathrel {\vbox{\hrule
   width\TriCDarrw@\vskip3\ex@\hrule width
   \TriCDarrw@}}}
\def\>##1>>{\setbox\z@\hbox{$\scriptstyle
 \;{##1}\;\;$}\global\bigaw@\TriCDarrw@
 \ifdim\wd\z@>\bigaw@\global\bigaw@\wd\z@\fi
 \hskip.5em
 \mathrel{\mathop{\hbox to \TriCDarrw@
{\rightarrowfill}}\limits^{##1}}
 \hskip.5em}
\def\<##1<<{\setbox\z@\hbox{$\scriptstyle
 \;{##1}\;\;$}\global\bigaw@\TriCDarrw@
 \ifdim\wd\z@>\bigaw@\global\bigaw@\wd\z@\fi
 \mathrel{\mathop{\hbox to\bigaw@{\leftarrowfill}}\limits^{##1}}
 }
 \CD@true\vcenter\bgroup\relax\let\\=\cr\iffalse}\fi
 \tabskip\z@skip\baselineskip20\ex@
 \ifTriV@
 \halign\bgroup
 &\hfill$\m@th##$\hfill\cr
#1&\multispan3\hfill$#2$\hfill&#3\\
&#4&#5\\
&&#6\cr\egroup%
\else
 \halign\bgroup
 &\hfill$\m@th##$\hfill\cr
&&#1\\%
&#2&#3\\
#4&\multispan3\hfill$#5$\hfill&#6\cr\egroup
\fi}
\def\@endTriCD{\egroup}
\newcommand{\sA}{{\mathcal A}}
\newcommand{\sB}{{\mathcal B}}
\newcommand{\sC}{{\mathcal C}}

\newcommand{\sE}{{\mathcal E}}
\newcommand{\sF}{{\mathcal F}}

\newcommand{\sH}{{\mathcal H}}

\newcommand{\sL}{{\mathcal L}}
\newcommand{\sM}{{\mathcal M}}
\newcommand{\sN}{{\mathcal N}}
\newcommand{\sO}{{\mathcal O}}

\newcommand{\sT}{{\mathcal T}}

\newcommand{\sX}{{\mathcal X}}

\newcommand{\A}{{\mathbb A}}

\newcommand{\C}{{\mathbb C}}
\renewcommand{\E}{{\mathbb E}}

\newcommand{\G}{{\mathbb G}}

\renewcommand{\L}{{\mathbb L}}
\newcommand{\M}{{\mathbb M}}

\newcommand{\BP}{{\mathbb P}}
\newcommand{\Q}{{\mathbb Q}}
\newcommand{\R}{{\mathbb R}}

\newcommand{\T}{{\mathbb T}}
\newcommand{\U}{{\mathbb U}}
\newcommand{\V}{{\mathbb V}}
\newcommand{\W}{{\mathbb W}}

\newcommand{\Z}{{\mathbb Z}}
\newcommand{\rk}{{\rm rank}}
\newcommand{\Hg}{{\rm Hg}}

\newcommand{\Sl}{{\rm Sl}}
\newcommand{\Gl}{{\rm Gl}}
\newcommand{\Sp}{{\rm Sp}}
\newcommand{\MT}{{\rm MT}}
\newcommand{\Mon}{{\rm Mon}}
\newcommand{\Gal}{{\rm Gal}}
\makeindex             

\begin{document}
\title*{Special families of curves, of Abelian varieties and of
certain minimal manifolds over curves}
\titlerunning{Families of manifolds over curves}
\author{Martin M\"oller\inst{1}, Eckart Viehweg\inst{1}\and Kang Zuo\inst{2}\thanks{This work has been supported by the ``DFG-Schwerpunktprogramm
Globale Methoden in der Komplexen Geometrie'', and by the DFG-Leibniz program.}}
\institute{Universit\"at Duisburg-Essen, Mathematik, 45117 Essen, Germany
  \texttt{martin.moeller@uni-essen.de}\\
  \texttt{viehweg@uni-essen.de}
\and Universit\"at Mainz,
Fachbereich 17, Mathematik,
55099 Mainz, Germany
\texttt{kzuo@mathematik.uni-mainz.de}}
\maketitle
\section*{Introduction}\label{in}

Let $f:X\to Y$ be a surjective morphism from an $n+1$ dimensional complex 
projective manifold $X$ to a curve $Y$. In this article we want to present some recent results about the structure of $f$, in particular about its discriminant locus $S$, i.e. the set of points $s$ in $Y$ with $f^{-1}(s)$ singular. So
$$
f:V=f^{-1}(U)\>>> U = Y\setminus S
$$
will always be smooth. We will usually require $\Delta=f^*S$ to be a normal crossing divisor.

If $Y=\BP^1$ and if $X$ is a curve of genus $g(X) > 0$ the Hurwitz formula implies that
$\#S \geq 3$. In the extremal case $\#S=3$, $X$ is obviously defined over a
number field. G.V. Bely{\u\i} has shown in \cite{B79} that
the existence of $f$ with $\#S=3$ characterizes curves defined over a number field. 

As conjectured by F. Catanese and M. Schneider, the first part can be generalized
for $n>0$. \vspace{.2cm}

\begin{theorem}[\cite{VZ01}]\label{in.1} Let $X$ be a complex projective manifold of Kodaira dimension $\kappa(X)\geq 0$ and let $f:X\to \BP^1$ be a surjective morphism. 
Then $f$ has at least $3$ singular fibres.
\end{theorem}
It is more reasonable to assume that the general fibre $F$ of $f$ is connected, or as we will say, that $f:X\to Y$ is
a family of $n$-dimensional manifolds, and to use the Hurwitz formula for the Stein factorization to get hold of the general case. $Y$ will be again a curve of 
arbitrary genus. Putting together results due to 
Parshin-Arakelov, Migliorini, Kov\'acs, Bedulev-Viehweg, Oguiso-Viehweg and Viehweg-Zuo
(see \cite{VZ01} and the references given there) one has:
\begin{theorem}\label{in.2} 
Let $f:X\to Y$ be a non-isotrivial family of $n$-folds, with general fibre $F$.
Assume either
\begin{enumerate}
\item[a.] $\kappa(F)=\dim(F)$, or
\item[b.] $\omega_F$ semiample.
\end{enumerate}
Then $\deg(\Omega^1_Y(\log S))=2g(Y)-2+\#S > 0$.
\end{theorem}
Here $f$ is isotrivial, if over some finite covering $Y'\to Y$
the pullback $X\times_YY'$ is birational to $F\times Y'$.

As a byproduct, the proof of \ref{in.2} in \cite{VZ01} (see also \cite{K02}) gives some explicit lower bounds for the degree of $\deg(\Omega^1_Y(\log S))$.  Writing $\delta$ for the number of singular
fibres of $f$, which are not semistable, i.e. not reduced normal crossing divisors, one finds
constants $\nu$ and $e$, depending only on the Hilbert polynomial $h$ of some polarization 
of the family, with
$$
\frac{\deg (f_* \omega^{\nu}_{X/Y})}{{\rm rank} (f_*
\omega^{\nu}_{X/Y})} \leq ( n \cdot (2 g(Y) - 2 + \#S) + \delta ) \cdot \nu
\cdot e.
$$
In Section \ref{co} we will use \cite[3.4]{VZ05} to get a similar bound
for arbitrary semistable families.

\begin{theorem} \label{in.3}
Assume that $f:X\to Y$ is a semistable family of $n$-folds. If $Y=\BP^1$ assume in addition that
$\#S\geq 2$. Then for all $\nu\geq 1$ with
$f_*\omega_{X/Y}^\nu\neq 0$
$$
\frac{\deg (f_* \omega^{\nu}_{X/Y})}{{\rm rank} (f_*
\omega^{\nu}_{X/Y})}\leq \frac{n\cdot \nu}{2} \cdot (2 g(Y) - 2 + \#S)=
\frac{n\cdot \nu}{2} \cdot \deg(\Omega^1_Y(\log S)).
$$
\end{theorem}
The assumption $\# S \geq 0$ is just added to guaranty that $\deg(\Omega^1_Y(\log S))$ is non-negative, and that the fundamental group of $Y\setminus S$ sufficiently large. Of course it will always hold true if one declares one or two smooth fibres to be singular.

There are families of Abelian varieties and of Calabi-Yau manifolds, with
$$
\deg (f_* \omega^{\nu}_{X/Y}) =\frac{\deg (f_* \omega^{\nu}_{X/Y})}{{\rm rank} (f_*
\omega^{\nu}_{X/Y})}= \frac{n\cdot \nu}{2} \cdot \deg(\Omega^1_Y(\log S)).
$$
We do not know, whether there are families of manifolds of higher Kodaira dimension with
such an equality.

Theorem \ref{in.3} implies Theorem \ref{in.2}. In fact, if $\deg(\Omega^1_Y(\log S))\leq 0$ the open set
$U$ is either an elliptic curve or it contains $\C^*$. In the first case the family is smooth, in the second one
one finds a new family over an \'etale covering $\C^* \to \C^*$, which is semistable. 
For a semistable family $f:X\to Y$ Theorem \ref{in.3} implies that $\deg (f_* \omega^{\nu}_{X/Y}) \leq 0$.
On the other hand, the assumptions a) or b) in Theorem \ref{in.2} imply that for $\nu \geq 2$ 
$$
\deg (f_* \omega^{\nu}_{X/Y}) \geq 0,
$$ 
and that $\deg (f_* \omega^{\nu}_{X/Y}) = 0$ if and only if $f$ is isotrivial. The same argument is used to show Theorem \ref{in.1}:

If there is a morphism with $2$ singular fibres, one uses the Hurwitz formula to show that the Stein factorization of $X\to \BP^1$ is $\BP^1$, and that the induced morphism has again at most $2$ singular fibres. Here
the positivity of $\deg (f_* \omega^{\nu}_{X/Y})$ follows from $\kappa(X)\geq 0$, and as above one obtains a contradiction
to the inequality stated in Theorem \ref{in.3}.

In addition, the bound in Theorem \ref{in.3} implies some weak analogue of the Shafarevich
conjecture on the finiteness of smooth families of curves over a fixed curve $U$. 

\begin{corollary}\label{in.4} Fix $(Y,S)$ and a polynomial $h$. Then up to deformation there are only finitely many non-isotrivial families $f:X\to Y$ smooth over $U=Y\setminus S$, such that:
\begin{enumerate}
\item $\omega_F$ is ample with Hilbert polynomial $h$.
\item $\omega_F^\rho=\sO_F$ for some $\rho>0$, and the smooth part $V\to U$
allows a polarization with Hilbert polynomial $h$.
\end{enumerate}
\end{corollary}
In fact, consider the moduli scheme $M_h$ parameterizing either canonically
polarized manifolds, or polarized manifolds with a torsion
canonical divisor, and with Hilbert polynomial $h$. In both cases
one can show (see \cite{V05}) that for a given $\nu\geq 2$ with $h(\nu)\neq 0$ there is a projective compactification
$\bar{M}_h$ and for some $p$ an invertible sheaf $\lambda_\nu^{(p)}$ on $\bar{M}_h$
with:
\begin{enumerate}
\item $\lambda_\nu^{(p)}$ is numerically effective, and ample with respect to $M_h$. 
\item If $f:X\to Y$ is a semistable family over a curve $Y$, whose general fibre belongs to the moduli functor, and if $\varphi:Y\to \bar{M}_h$ is the extension of the induced morphism
to $M_h$, then 
$$
\varphi^*\lambda_\nu^{(p)}=\det(f_* \omega^{\nu}_{X/Y})^p.
$$
\end{enumerate}
Here an invertible sheaf $\lambda$ is ``ample with respect to $M_h$'' if there exists some modification $\tau:\bar{M}'_h\to \bar{M}_h$, some $\mu\gg 1$ and an effective divisor $E$ on $\bar{M}'_h$, supported in $\tau^{-1}(\bar{M}_h\setminus M_h)$, such that $\tau^*\lambda\otimes\sO_{\bar{M}'_h}(-E)$ is ample.

So Theorem \ref{in.3} says that for given $Y$ and $S$ the degree of
$\varphi^*\lambda_\nu^{(p)}$ is bounded, which implies Corollary \ref{in.4}.
Using a different argument (see \cite[6.2]{VZ02}) one can extend \ref{in.4}
to families of minimal models of arbitrary (non-negative) Kodaira dimension.

As we will see in Section \ref{co} the Theorems \ref{in.1}, \ref{in.2} and
\ref{in.3} follow from the study of variations of Hodge structures of certain families
obtained as cyclic coverings of $X\to Y$.\vspace{.2cm}

It remains the question, what is special about families with few singular fibres.
Here the only results are for families controlled by their variation
of Hodge structures $R^nf_*\C_V$ (or $R^1f_*\C_V$ for families of Abelian varieties).
Moreover, we will assume that the family has semistable reduction, or slightly
weaker, that the local monodromies in all $s\in S$ are unipotent. Before
stating the results, and before giving ``few singular fibres'' a precise meaning, we will need some notations.

Let $\V\subset R^k\C_V$ be a $\C$ sub variation of Hodge structures.
$\V\otimes_C \sO_U$ extends to a locally free sheaf $\sH$ on $Y$ in such a way that the Gau{\ss}-Manin connection acquires logarithmic singularities. We choose for $\sH$ the Deligne extension, i.e. an extension such that the real part of the local residues are zero.
The Hodge filtration extends to a holomorphic filtration on $\sH$, and the extended Gau{\ss}-Manin connection defines on the associated graded bundle the structure of a logarithmic Higgs bundle $(F,\tau)$, i.e. a locally free sheaf $F$ together with a collection of maps
$$
\tau_{p,q}:F^{p,q}\>>> F^{p-1,q+1}\otimes \Omega_Y^1(\log S), \ \ \ \ \ \ p+q=k.
$$
The bundle maps $\tau_{p,q}$ can be iterated to obtain
$$
\tau^{(\ell)}: F^{k,0} \>>> F^{k-1,1}\otimes \Omega^1_Y(\log S)\>>> \cdots
\>>> F^{k-\ell,\ell}\otimes S^\ell(\Omega^1_Y(\log S)).
$$
\begin{definition}\label{in.5} \ 
\begin{enumerate}
\item[i.] We call $\tau^{(k)}:F^{k,0} \to F^{0,k}\otimes S^k(\Omega^1_Y(\log S))$
the Griffiths-Yukawa coupling of $\V$ (or of $f$ in case $\V=R^kf_*\C_V$).
\item[ii.] The Higgs field is strictly maximal, if $F^{k,0}\neq 0$ and if the $\tau_{p,q}$
are all isomorphisms.
\end{enumerate}
\end{definition}
As it will turn out, the property ii) is of numerical nature.
\begin{lemma}[\cite{VZ05}]\label{in.6}
Assume that $\V$ is a variation of polarized complex Hodge structures 
of weight $n$ with unipotent local monodromy in all $s\in S$, and with
logarithmic Higgs bundle $(\bigoplus F^{p,q},\tau_{p,q})$. If $F^{k,0}\neq 0$ one has the
Arakelov inequality 
\begin{equation}\label{eqin.1}
\deg(F^{k,0}) \leq \frac{k}{2}\cdot\rk(F^{k,0})\cdot\deg(\Omega^1_Y(\log S)),
\end{equation}
and (\ref{eqin.1}) is an equality if and only if one has a decomposition $\V=\V_1\oplus \V_2$ where the Higgs field of $\V_1$ is strictly maximal and where
$\V_2$ is a variation of polarized complex Hodge structures, zero in bidegree
$(k,0)$. 
\end{lemma}
For families of curves, the inequality (\ref{eqin.1}) is due to Arakelov, and for families of Abelian varieties it has been shown by Faltings in \cite{F83}.

If $\#S$ is even, \cite[3.4]{VZ03} gives a more precise description of
$\V_1$. Choose a logarithmic theta characteristic, 
i.e. an invertible sheaf $\sL$ with $\sL^2=\Omega_Y^1(\log S)$. The existence
of $\V$ in \ref{in.6} implies that $\deg(\Omega^1_Y(\log S)>0$. 
Then the Higgs bundle $\sL\oplus \sL^{-1}$ with Higgs field
$$
\sL \>\tau> \simeq > \sL^{-1}\otimes\Omega_Y^1(\log S) 
$$
is stable. As recalled in Section \ref{hb} it must be the Higgs bundle of a polarized variation of Hodge structures $\L$ of weight one.
$\L$ is unique up to the tensor product with a unitary rank one local system,
induced by a two-division point of ${\rm Pic}^0(Y)$.
\begin{lemma}\label{in.7}
Assume that $\#S$ is even, and that $\V_1$ is a variation of Hodge structures with a strictly maximal Higgs field and with unipotent local monodromies. Then there exists
a unitary local system $\T$ on $Y$, regarded as a variation of Hodge structures in bidegree $(0,0)$, with $\V_1=S^k(\L)\otimes_C \T|_U$.
\end{lemma}

So ``special families'' or ``families with few singular fibres'' will be those
where the variation of Hodge structures $R^k\C_V$ contains an irreducible sub variation $\V$ with a strictly maximal Higgs field, or equivalently a sub variation $\V$ for whose Higgs bundle the equality
\begin{equation}\label{eqin.2}
\deg(F^{k,0}) = \frac{k}{2}\cdot\rk(F^{k,0})\cdot\deg(\Omega^1_Y(\log S))
\end{equation}
holds. Here $k=n$, except for families of Abelian varieties, where we usually choose $k=1$.\vspace{.2cm}

Let us first consider a semistable family $f:X\to Y$ of Abelian varieties, smooth over $U=Y\setminus S$,
and $V=f^{-1}(U)$.
\begin{theorem}[\cite{VZ04}]\label{in.8}
Assume that each irreducible and non-unitary sub variation $\V$ of Hodge structures in $R^1f_*\C_V$
has a strictly maximal Higgs field. Then (replacing $U$ by an \'etale covering and $V$ by the pullback family) $U$ is a Shimura curve of Hodge type, and $f:V\to U$ the corresponding universal family.
\end{theorem}   
The construction of Shimura curves will be recalled in Section \ref{mh}. Let us just mention here, that a Shimura curve $U$ is an \'etale covering of a certain moduli space of Abelian varieties with prescribed Mumford-Tate group
and a suitable level structure. So whenever we talk about Shimura curves, this is a property up to \'etale coverings of $U$. This allows in particular to assume that $\# S$ is even.

As we will recall in Section \ref{tc} the converse of Theorem \ref{in.8} was shown in \cite{Moe3}:
\begin{theorem}\label{in.9}
If $f:V\to U$ is the universal family over a Shimura curve its Higgs field is strictly maximal.
\end{theorem}

Because of Koll\'ar's decomposition \cite{K87} (see Lemma \ref{hb.9} below) one can write $f_*\Omega_{X/Y}^1(\log S)$ as a direct sum
of an ample sheaf $\sA$ and a subsheaf $\sB$, flat for the Gau{\ss}-Manin connection.
Correspondingly on the local system side one has a decomposition 
$R^1f_*\C_V=\W\oplus\U$, where $\U$ is unitary and invariant under complex conjugation, and where the $(1,0)$ component of $\W$ is ample. Hence the Higgs bundle corresponding to $\W$ is of the form
$(\sA\oplus \sA^\vee,\theta)$ and using \ref{in.6} one can restate the Theorem \ref{in.8} as:\\[.2cm] 
{\it If $\deg(\sA) = \frac{1}{2}\cdot\rk(\sA)\cdot\deg(\Omega^1_Y(\log S))$, then
$U$ is a Shimura curve of Hodge type, and $f:V\to U$ the corresponding universal family.}\\[.2cm]
The Lemma \ref{in.7} implies that the Arakelov equality forces $\W=\L\otimes_\C \T$ for a unitary system $\T$.
In particular, if $S\neq \emptyset$, replacing $Y$ by an \'etale covering, one can assume that
$\T=\C_U^{\oplus r}$ and $\U=\C_U^{\oplus n-r}$. As in \cite{VZ04} one deduces the next Corollary.
\begin{corollary}\label{in.10}
If $S\neq \emptyset$ consists of an even number of points, and if  
$$
\deg(\sA) = \frac{1}{2}\cdot\rk(\sA)\cdot\deg(\Omega^1_Y(\log S)),
$$
then (replacing $Y$ by an \'etale covering) $f:X\to Y$ is isogenous to $$E\times_Y\cdots \times E_Y \times B,$$
where $B$ is an Abelian variety of dimension $n-r$ and where $E\to Y$ is a modular family of elliptic curves.
\end{corollary}
The proof of Corollary \ref{in.10} is easy if one assumes that the maximal unitary subsystem
of $R^1f_*\C_V$ is defined over $\Q$, and it will be sketched in Section \ref{mh}. 
This condition holds true if $S\neq\emptyset$. 
For $S=\emptyset$ the proof of Theorem \ref{in.8} in \cite{VZ04} also gives back the known classification of Shimura curves. 

As a supplement to \ref{in.8} we will show in Section \ref{mh}:
\begin{corollary}\label{in.11}
If the maximal unitary local subsystem $\U$ of $R^1f_*\C_V$ is trivial, then
the family $f:V\to U$ is rigid. In particular it is defined over $\bar{\Q}$.
\end{corollary}
Let us consider next families of curves of genus $g\geq 2$, and let us return to the
local systems $\L$, defined by logarithmic theta characteristics $\sL$.
Obviously $\L$ has a strictly maximal Higgs field.
\begin{theorem}[\cite{Moe2}]\label{in.12}
Let $f:X\to Y$ be a semistable family of curves of genus $g\geq 2$, smooth over $U=Y\setminus S$, and $V=f^{-1}(U)$.
\begin{enumerate}
\item $R^1f_*\C_V$ contains a rank two sub variation of Hodge structures $\L$ with a strictly maximal Higgs field
if and only if $U$ is a Teichm\"uller curve, and $V\to U$ the corresponding universal family. 
\item Teichm\"uller curves are defined over number fields.
\end{enumerate}
\end{theorem}
The definition of a Teichm\"uller curve will be given in Section \ref{tc}.
Roughly speaking, one considers geodesics in the Teichm\"uller space,
constructed by an ${\rm Sl}(2,\R)$-action on the real and imaginary part of a given holomorphic differential form. If the quotient by a suitable lattice in ${\rm Sl}(2,\R)$ is an algebraic curve, it is called a Teichm\"uller curve in $M_g$. 
As for Shimura curves, we always allow ourselves to replace the lattice by a smaller one, hence
the Teichm\"uller curve by an \'etale cover. We say that $f:V\to U$ is the universal family 
if the morphism $U\to M_g$ is induced by $f$.

It is a striking fact, that in spite of the differential geometric definition of Teichm\"uller curves, they can be characterized analytically, or even by a numerical property of the Higgs bundle. 
This has a number of applications in the theory of Teichm\"uller curves (see \cite{Moe2} and \cite{Moe3}), and it allows to construct new examples of such special curves in $M_g$.

At the same time, Theorem \ref{in.12} allows to use differential geometric properties
of Teichm\"uller curves to study variations of Hodge structures for families
of curves. For example, at the present moment the only known proof of the next Corollary  
relies on Theorem \ref{in.12}, and on differential geometric methods.
\begin{corollary}\label{in.13} In Theorem \ref{in.12} the variation of Hodge structures
$R^1f_*\C_V$ can not contain two different sub variation of Hodge structures $\L_1$ and $\L_2$, both with a strictly maximal Higgs field.  
\end{corollary}
There remains the question whether there are higher rank irreducible sub variations of Hodge structures $\V$ with a strictly maximal Higgs field. By Lemma \ref{in.7} those are of the form $\L\otimes_\C \T$, where $\T$ is an irreducible unitary system, at least after replacing $U$ by an \'etale covering of degree two.
As we will see in Section \ref{hb} $\rk(\T)\geq 2$ implies that $S=\emptyset$, and squeezing the methods used
to prove \ref{in.2} a bit more, one can show that this is not possible (see \cite{VZ05}). So putting everything together one obtains:
\begin{corollary}\label{in.14}
Let $f:X\to Y$ be a semistable family of curves of genus $g\geq 2$, smooth over $U=Y\setminus S$.
Then $R^1f_*\C_V$ does not contain a sub variation of Hodge structures $\V$ with a strictly maximal Higgs field and with $\rk(\V)>2$. 

If $S=\emptyset$ then $R^1f_*\C_V$ does not contain any sub variation of Hodge structures $\V$ with a strictly maximal Higgs field.
\end{corollary}
As we will see in section \ref{tc}, the rank two local subsystem $\L\subset R^1f_*\C_V$ can be defined over a number field $K$, i.e. it is of the form $\L_K\otimes_K\C$ for a local subsystem $\L_K\in R^1f_*K_V$, but in general
it will not be defined over $\Q$. Replacing $K$ by its Galois hull, for $\sigma\in {\rm Gal}(K/\Q)$ 
the conjugate local system $\L^\sigma$ can only have a strictly maximal Higgs field, if
$\L=\L^\sigma$. The defect in the inequality (\ref{eqin.1}) for $\L^\sigma$
has been studied in \cite{BM05}.\vspace{.2cm}

In general things are changing, if one replaces $f:X\to Y$ by its family of Jacobians. Passing from a family of curves to its Jacobian, one might replace the discriminant locus $S$ by a subset $S'$. Singular fibres with a compact Jacobian, for example those with two components meeting in one point, are not contributing to the discriminant locus of the family of Jacobians. If $f:V\to U$ is the universal family over a Teichm\"uller curve, then the existence of
the local subsystem $\L$ together with the inequality (\ref{eqin.1}) applied to the family of Jacobians,
shows that $S'=S$, hence that all singular fibres of $f$ have a non-compact Jacobian. One can ask, whether
there are any Teichm\"uller curves $U\subset M_g$ in the moduli space of curves of genus $g$ which are Shimura curves in the moduli space $A_g$ of Abelian varieties (with a suitable level structure). By \cite{Moe3} the answer is ``no'' for
$g\geq 4$  and for $g=2$, whereas for $g=3$ there exists one example, essentially unique.

On the other hand, by the Corollaries \ref{in.13} and \ref{in.14} we know, that
the only curves $U\subset M_g$, which are Shimura curves in $A_g$ have to be induced by
a family $f:X\to Y$ with $R^1f_*\C_V=\L\oplus \U$ with $\rk(\L)=2$ and with $\U$ unitary. 
This implies that $U$ is a Teichm\"uller curve, and that $U$ is not projective.

\begin{theorem}\label{in.15} \
\begin{enumerate}
\item[a.] For $g\geq 2$ the moduli space $M_g$ does not contain any compact Shimura curve.
\item[b.] $M_g$ contains a non-compact Shimura curve $U$ if and only if $g=3$. 
The Shimura curve $U\in M_3$ is a Teichm\"uller curve, essentially unique.
\end{enumerate}
\end{theorem}
In this article, we will explain the main ingredients used in the introduction, and we will sketch the proofs of some of the results mentioned above. In Section \ref{hb} we explain the Simpson correspondence for
variations of Hodge structures over curves, and we will prove Lemma \ref{in.6}
and Lemma \ref{in.7} for $k=1$. For $k>1$, Lemma \ref{in.6} will be shown under the additional assumption that $\rk(F^{k,0})=1$.

In Section \ref{co} we prove Theorem \ref{in.3}, and we give some hints, how one 
obtains Corollary \ref{in.14} from Corollary \ref{in.13}. 

The characterization of Shimura curves will be discussed in Section \ref{mh}.
We prove Theorem \ref{in.8} for families where 
$R^1f_*\C_V$ itself has a strictly maximal Higgs field, hence we exclude
unitary direct factors. We sketch the proofs of Corollaries \ref{in.10} and \ref{in.11}
and we also reproduce the proof of Theorem \ref{in.9}
from \cite{Moe3}.

Section \ref{tc} gives an introduction to the theory of Teichm\"uller curves.
In particular we sketch the proof of Theorem \ref{in.12} 
and of the Corollary \ref{in.13}.

Finally the non-existence of Teichm\"uller curves in $M_g$, for $g=2$ and $g\geq 4$
and the proof of Theorem \ref{in.15}, b), will be discussed in Section \ref{to}.
\section{Higgs bundles over curves and Arakelov inequalities}\label{hb}
We will frequently use C. Simpson's correspondence between
polystable logarithmic Higgs bundles of degree zero and representations of
the fundamental group $\pi_1(U,*)$. Recall that a logarithmic Higgs bundle
is a locally free sheaf $E$ on $Y$ together with an $\sO_Y$ linear morphism
$\theta:E\to E\otimes \Omega^1_Y(\log S)$ with $\theta\wedge\theta=0$. The usual definitions of stability (and semistability) for locally free sheaves extend to Higgs bundles,
by requiring that
$$
\mu(F)=\frac{\deg(F)}{\rk(F)} < \mu(E)=\frac{\deg(E)}{\rk(E)} 
$$
(or $\mu(F)\leq \mu(E)$) for all subsheaves $F$ with $\theta(F)\subset F\otimes \Omega^1_Y(\log S)$.
\begin{theorem}[C. Simpson \cite{S90}]\label{hb.1}
There exists a natural equivalence between the category of direct
sums of stable filtered regular Higgs bundles of degree zero, and
of direct sums of stable filtered local systems of degree zero.
\end{theorem}
We will not recall the definition of a ``filtered regular'' Higgs bundle
\cite[page 717]{S90}, and just remark that for a Higgs bundle corresponding
to a local system $\V$ with unipotent monodromy around the points in
$S$ the filtration is trivial, and automatically $\deg(\V)=0$.

In general it is impossible to describe the Higgs bundle $(E,\theta)$
explicitly in terms of the corresponding local system $\V$, with two exceptions: 
\begin{example}\label{hb.2}
If $\T$ is a unitary local system on $U$ with unipotent monodromy operators in $s\in S$,
then the corresponding Higgs bundle is the Deligne extension $N$ of $N_0=\T\otimes \sO_U$
and the Higgs field $\theta$ is the zero map. 
\end{example}
The second example already occurred in the Introduction.
\begin{example}\label{hb.3}
Let $\V$ be a polarized $\C$ variation of Hodge structures of weight $k$ and with
unipotent local monodromy operators. The $\sF$-filtration on $F_0=\V\otimes_\C \sO_U$
extends to a locally splitting filtration on the Deligne extension $F$
$$
\sF^{k+1} \subset \sF^{k} \subset \cdots \subset \sF^0.
$$
We will usually assume that $\sF^{k+1}=0$ and $\sF^0=F$.
The Griffiths transversality condition for the Gau{\ss}-Manin connection
$\nabla$ says that
$$
\nabla(\sF^p)\subset \sF^{p-1}\otimes \Omega^1_Y(\log S),
$$
and hence $\nabla$ induces a $\sO_Y$ linear map
$$
\theta_{p,k-p}: F^{p,k-p}=\sF^p/\sF^{p+1} \>>> F^{p-1,k-p+1}=\sF^{p-1}/\sF^{p}\otimes
\Omega^1_T(\log S).
$$
So
$$
\big(F=\bigoplus_p F^{p,k-p}, \theta=\bigoplus \theta_{p,k-p}\big)
$$
is a Higgs bundle, and it is the image of $\V$ under the Simpson correspondence
in Theorem \ref{hb.1}.
\end{example}
\begin{properties}\label{hb.4}
Let $\V$ be a direct sum of irreducible $\C$ local system and let
$(F,\theta)$ be the corresponding Higgs bundle. 
\begin{enumerate}
\item $\displaystyle F=\bigoplus_{p+q=k}F^{p,q}$ with $\theta(F^{p,q})\subset
F^{p-1,q+1}\otimes\Omega^1_Y(\log S)$ if and only if $\V$ is a local system underlying
a complex polarized variation of Hodge structures.
\item Under the equivalent conditions in 1) let $\T$ be a local subsystem of $\V$ with Higgs bundle
$$
\big(N=\bigoplus_p N^{p,k-p},0\big) \ \subset \ \big(F=\bigoplus_p F^{p,k-p}, \theta=\bigoplus \theta_{p,k-p}\big).
$$
Then $\T$ is unitary.
\end{enumerate}
\end{properties}
\begin{proof}
The condition on $(F,\theta)$ in 1) is just saying, that it is a system of Hodge bundles,
as defined in \cite{S88}. So the first part follows from \cite{S88} and \cite{S90}.

Let $\Theta(N,h)$  denote the curvature of
the Hodge metric $h$ on $F$ restricted to $N,$
then by \cite{G84}, chapter II we have
$$
\Theta(N,h|_N)=-\theta_N\wedge\bar\theta_N-\bar\theta_N\wedge\theta_N=0.
$$
This means that $h|_N$ is a flat metric. Hence,  $\T$ is a
unitary local system.
\myqed  \end{proof}
The Simpson correspondence does not say anything about the field of definition for $\V$.
Here we say that $\V$ is defined over a subfield $K$ of $\C$ if there is a $K$-local system
$\V_K$ with $\V=\V_K\otimes_K\C$.
In different terms, for $\mu=\rk(\V)$ the representation
$$
\gamma_\V:\pi_0(U,*) \>>> {\rm Gl}(\mu,\C)
$$
is conjugate to one factoring like
$$
\gamma_\V:\pi_0(U,*) \>>> {\rm Gl}(\mu,K) \>>> {\rm Gl}(\mu,\C).
$$
If $\V$ is defined over $K$, and if $\sigma:K\to K'$ is an
isomorphism, we will write $\V_K^\sigma$ for the local system
defined by
$$
\gamma_\V:\pi_0(U,*) \>>> {\rm Gl}(\mu,K) \> \sigma >> {\rm
Gl}(\mu,K'),
$$
and $\V^\sigma=\V_K^\sigma\otimes_{K'} \C$.

By \cite{D71} (for $K=\Q$) and by \cite{D87} the category of polarized $K$-variations of Hodge structures is semisimple. By \cite[Proposition 1.13]{D87} one has:
\begin{lemma}\label{hb.5}
A local system $\V$, underlying a polarized variation of Hodge structures,
decomposes as 
$$
\V = \bigoplus_{i=1}^r (\V_i \otimes W_i),$$
where $\V_i$ are pairwise non-isomorphic 
irreducible $\C$-local systems and $W_i$
are non-zero $\C$-vector spaces.

Moreover the $\V_i$ and the  $W_i$ carry 
polarized variations of Hodge structures, whose tensor product and sum 
gives back the Hodge structure 
on $\V$. The Hodge structure on the $\V_i$ (and $W_i$) 
is unique up to a shift of the bigrading.
\end{lemma}
Suppose that $\W$ is a local system defined over a number field $L$. The local system $\W_{L}$
is given by a representation $\rho: \pi_1(U,*) \to
\text{Gl}(W_{L})$ for the fibre $W_{L}$ of $\W_{L}$ over the base
point $*$.

Fixing a positive integer $r<\mu$ let $\mathcal G(r,\W)$ denote the
set of all rank-r local subsystems of $\W$ and let ${\rm
Grass}(r,W_L)$ be the Grassmann variety of $r$-dimensional subspaces.
Then $\mathcal G(r,\W)$ is the subvariety of
$$
{\rm Grass}(r, W_{L})\times_{{\rm Spec}(L)}{\rm Spec}(\C)
$$
consisting of the $\pi_1(U,*)$-invariant points. In particular,
it is a projective variety defined over $L$. A $K$-valued point
of $\mathcal G(r,\W)$ corresponds to a local subsystem of
$\W_K=\W_{L}\otimes_L K$. One obtains the following well known property.

\begin{lemma}\label{hb.6}
If $[\V]\in \mathcal G(r,\W)$ is an isolated point,
then $\V$ is defined over $\bar\Q.$
\end{lemma}

\begin{lemma}\label{hb.7}
Let $\W$ be a polarized variation of Hodge structures
defined over $L$, and let $\V\subset \W$ be an irreducible local subsystem
of rank $r$ defined over $\C,$. Then $\V$ can be deformed to a local subsystem
$\V_t\subset\W$, which is isomorphic to $\V$ and which is defined over
a finite extension of $L$.
\end{lemma}
\begin{proof} By Lemma \ref{hb.5} $\W$ is completely reducible over $\C$. Hence
we have a decomposition $\W=\V\oplus\V'$.

The space $\mathcal G(r,\W)$ of rank $r$ local subsystems of $\W$
is defined over $L$ and the subset
$$ \{\V_t\in \mathcal G(r,\W); \mbox{ the composite }
\V_t \subset \V\oplus \V' \> pr_1 >> \V\mbox{ is non zero }\}$$
forms a Zariski open subset. So there exists some $\V_t$ in this subset,
which is defined over some finite extension of $L$. Since $p: \V_t\to \V$
is non zero, since ${\rm rank}(\V_t)= {\rm rank}(\V)$, and since
$\V$ is irreducible, $p$ is an isomorphism.
\myqed  \end{proof}

The Lemmata \ref{hb.6} and \ref{hb.7} are just the starting points to show,
that certain local subsystems of $R^kf_*\C_V$ are defined over number fields
(see \cite{VZ04}). Let us give a typical example:
\begin{lemma}\label{hb.8} Let $\W$ be defined over a (real) number field, and let $\W=\V\oplus\U$ be a decomposition such that the
Higgs field of $\V$ is maximal, and such that $\U$ is unitary. 
\begin{enumerate}
\item  Then $\V$ and $\U$ are defined over a (real) number field,
as well as the decomposition.
\item If $S\neq \emptyset$, then $\U$ is defined over $\Q$ and it trivializes over
an \'etale covering of $Y$.
\end{enumerate}
\end{lemma}
\begin{proof}
Consider a family $\V_t$ of local subsystems of $\W$ for $t$ in a small disk $\Delta$, and with $\V_0=\V$. Since
the Higgs field of $\V_0$ is maximal, one may assume that the one of $\V_t$ is maximal for all $t\in \Delta$.
Then the projection $\V_t\to \U$ must be zero. Otherwise, the complete reducibility of local systems
coming from variations of Hodge structures in \ref{hb.5} implies that $\V_t$ and $\U$ have a common factor,
obviously a contradiction. 

So $\V$ is rigid as a local subsystem, hence defined over $\bar{\Q}$. The same argument works for $\U$ instead of $\V$.

If $\W$ is defined over $\R$, then $\bar{\V}\to \U$ again has to be the zero map. So $\V$ and $\U$ have to be defined over $\R$.

Assume now that $S\neq\emptyset$. By (1) we know that the decomposition is defined over a number field. For the local monodromy operators ``unitary and unipotent'' implies that the 
nilpotent part of the monodromy is zero. This is invariant under conjugation by $\sigma\in {\rm Gal}(\bar{\Q}/\Q)$, and $\U^\sigma=\U$. 

Finally, since $\U$ is a local subsystem of a variation of Hodge structures,
it carries a $\Z$-structure. Then the image of the corresponding monodromy representation is finite.\myqed
\end{proof}
Lemma \ref{hb.5} together with Theorem \ref{hb.1} implies that the Higgs
bundle $(F,\theta)$ of a variation of Hodge structures $\V$ is 
polystable (as a Higgs bundle) and of degree zero. As a first application one obtains Koll\'ar's decomposition of the
sheaf $F^{k,0}$, mentioned in the introduction.
\begin{lemma}\label{hb.9}
Let $(F=F^{k,0}\oplus \cdots \oplus F^{0,k},\theta)$ be the logarithmic Higgs bundle
of a variation of Hodge structures $\V$ with unipotent local monodromies.
Then $F^{k,0}=\sA\oplus \sB$ with $\sA$ ample and with $\sB$ flat.
\end{lemma}
\begin{proof}
Replacing $\V$ by one of its irreducible direct factors, one has to show that
$F^{k,0}$ is ample, provided that $\theta\neq 0$.
 
A quotient sheaf $\sN$ of $F^{k,0}$ gives rise to a quotient Higgs bundle
$(\sN,0)$ of $(F,\theta)$, hence \ref{hb.1} implies that $\deg(\sN)\geq 0$. 
If $\deg(\sN)= 0$, then it corresponds to a local subsystem, which is excluded by the irreducibility of $\V$. Hence all quotients of $F^{k,0}$ have strictly positive degree, which implies that $F^{k,0}$ is ample.
\myqed  \end{proof}
Let us consider for a moment the case $k=1$ and $\V=R^1f_*\C_V$. Then $(\sB\oplus \sB^\vee,0)$
is a sub Higgs bundle, hence by Lemma \ref{hb.4}, 2), it corresponds to a unitary subbundle of $\V$. 
The latter is defined over $\Q$ and by Lemma \ref{hb.8}, (1), this decomposition is
defined over $\bar\Q\cap \R$.

The Arakelov inequality stated in Lemma \ref{in.6} follows from the polystability
of the Higgs bundles (see \cite{VZ03} and \cite{VZ05}).
We will indicate the proof just in two simple cases, the one of weight one 
variations of Hodge structures, and the one where the $(k,0)$ part is one dimensional.\par\medskip \noindent
{\it Proof of Lemma \ref{in.6} for weight one.} \ 
Let $\sA\subset F^{1,0}$ be a subsheaf, and let $\sC\otimes\Omega_Y^1(\log S)$ be
its image under $\theta_{1,0}$. Then $\sA\oplus\sC$ is a Higgs
subbundle of $F^{1,0}\oplus F^{0,1}$, and \ref{hb.1} implies that
$\deg(\sA)+\deg(\sC)\leq 0$. Since $({\rm ker}(\theta|_\sA),0)$ is a sub Higgs bundle,
one has
\begin{multline}\label{eqhb.1}
\deg(\sA) \leq \deg(\sC) + \rk(\sC)\cdot\deg(\Omega_Y^1(\log S))\\
\leq \deg(\sC) + \rk(\sA)\cdot\deg(\Omega_Y^1(\log S))\\
\leq
-\deg(\sA) + \rk(\sA)\cdot\deg(\Omega_Y^1(\log S)),
\end{multline}
and
$$
\frac{\deg(\sA)}{\rk(\sA)} \leq \frac{1}{2}\deg(\Omega_Y^1(\log
S)).
$$
For $\sA=F^{1,0}$ one obtains the inequality (\ref{eqin.1}) in \ref{in.6}.

If (\ref{eqin.1}) is an equality, all the inequalities in (\ref{eqhb.1})
are equalities. This implies that $\rk(\sA)=\rk(\sC)$, hence that $\theta|_\sA$
is injective. Furthermore $\deg(\sA)+\deg(\sC)= 0$, hence $\sA\oplus \sC$ corresponds to a sub-local system of $\V$ with a maximal Higgs field.
\myqed \vspace{.3cm}
 
In fact we have shown a bit more. Assume that (\ref{eqin.1}) in Lemma \ref{in.6} is an equality.
If $\V$ is irreducible 
$$
\frac{\deg(\sA)}{\rk(\sA)} < \frac{1}{2}\deg(\Omega_Y^1(\log
S))= \frac{\deg(F^{1,0})}{g_0},
$$
except for $\sA=F^{1,0}$. Hence $F^{1,0}$ is stable. 
By duality one obtains the same for $F^{0,1}$. 

If $\V$ is not irreducible, applying this argument to all the direct factors of $\V$ one
finds:
\begin{lemma}\label{hb.10}
If $(F^{1,0}\oplus F^{0,1},\theta)$ is the logarithmic Higgs bundle
of a variation of Hodge structures $\V$ with unipotent local monodromies
and with a strictly maximal Higgs field, then the sheaves
$F^{1,0}$ and $F^{0,1}$ are polystable and the Higgs field $\theta:F^{1,0}\to F^{0,1}\otimes \Omega^1_Y(\log S)$ is a morphism between polystable sheaves of the same slope.
\end{lemma}
Although in this survey we only consider families over curves, let us
make one remark about families over a higher dimensional basis.
One can hope that the polystability of the Hodge bundles 
$F^{1,0}$ and $F^{0,1}$ is again enforced by a numerical condition on the degrees
with respect to $\omega_Y(S)$, and that it forces $U=Y\setminus S$ to be a Shimura
variety. Some first results in this direction have been obtained in \cite{VZ05b}. \par\medskip
\noindent{\it Proof of Lemma \ref{in.7} for $k=1$.} \ 
By assumption $\#S$ is even, so we can choose a logarithmic theta characteristic $\sL$
and an isomorphism 
$$
\tau:\sL\to \sL^{-1}\otimes \Omega^1_Y(\log S),
$$
giving us the $\C$-variation of Hodge structures $\L$.
By Lemma \ref{hb.10}, the sheaf $\sT=F^{1,0}\otimes \sL^{-1}$ is polystable of degree zero, hence
isomorphic to $\T\otimes \sO_Y$ for a unitary local system $\T$.
The isomorphism
$$
\sT\otimes \sL = F^{1,0}\> \theta > \simeq > F^{0,1}\otimes
\Omega^1_Y(\log S)\> >\simeq > F^{0,1}\otimes \sL^{2}
$$
induces an isomorphism $\phi: \sT\otimes \sL^{-1}\> \simeq >> F^{0,1}$,
such that $\theta=\phi\circ({\rm id}_\sT \otimes \tau)$. Hence
the Higgs bundles $(F^{1,0}\oplus F^{0,1},\theta)$ and
$(\sT\otimes(\sL\oplus \sL^{-1}),{\rm id}_\sT\otimes{\tau})$ are
isomorphic, and $\V\simeq \T\otimes_\C \L$.
\myqed
\par\medskip \noindent
{\it Proof of Lemma \ref{in.6} under the additional assumption $\rk(F^{k,0})=1$.} \ 
For later use, let us consider a slightly more general situation:\\[.1cm]
$\V$ is a variation of polarized complex Hodge structures and with unipotent local monodromy
in $s\in S$, and with logarithmic Higgs bundle $(\bigoplus F^{k,0}, \theta)$. Let $\sH \subset F^{k,0}$ be an invertible subsheaf, and let $E\subset F$ be the sub Higgs bundle generated by $\sH$.
Writing $T_Y(-\log S)$ for the dual of $\Omega^1_Y(\log S)$, for some $q_0$ 
$$
E^{k-q,q}=F^{k-q,q}\cap E=
\left\{ \begin{array}{ll}
\sH \otimes T_Y(-\log S)^{q} & \mbox{for \ } q \leq q_0\\
0 & \mbox{for \ } q >q_0.
\end{array}\right. .
$$
Theorem \ref{hb.1} implies that
$$
\deg(E)=(q_0+1)\cdot \deg(\sH) - \frac{q_0\cdot(q_0+1)}{2}\deg(\Omega^1_Y(\log S))
\leq 0,
$$
hence 
\begin{equation}\label{eqhb.4}
\deg(\sH)\leq\frac{q_0}{2}\deg(\Omega^1_Y(\log S))\leq\frac{k}{2}\deg(\Omega^1_Y(\log S)).
\end{equation}
In Lemma \ref{in.6} we choose $\sH=F^{k,0}$, and we obtain the Arakelov inequality. If this is an equality, then $q_0=k$ and $E$ corresponds to a local subsystem, by construction with a strictly maximal Higgs field.
\myqed   
\section{Coverings and bounds for subsheaves of the direct images}\label{co}
Let us recall from \cite{VZ05} the proof of the following proposition, as a first application of the theory of Higgs bundles.
\begin{proposition}\label{co.1}
Let $f:X\to Y$ be a semistable family of $n$-folds over a curve $Y$,and smooth over $U=Y\setminus S$.
If $Y=\BP^1$ assume that $\#S\geq 2$, and consider for $\nu \geq 1$ an invertible subsheaf $\sH$ of $f_*\omega_{X/Y}^\nu$. Then  
$$
\deg{\sH} \leq \frac{n\cdot \nu}{2}\cdot \deg(\Omega^1_Y(\log S)).
$$
\end{proposition}
Before sketching the proof of the proposition, let us show that it implies Theorem \ref{in.3}.
\par\medskip \noindent
{\it Proof of Theorem \ref{in.3}.} \
Since either $g(Y)\geq 1$, or $Y=\BP^1$ and $\# S \geq 2$, one can always find a covering $Y'\to Y$, which is  unramified over $U$ and which has a prescribed ramification order in $s\in S$.
Since $f$ is semistable, the inequality in Theorem \ref{in.3} is compatible with replacing $Y$ by $Y'$ and $f$ by a desingularization of the pullback family.

For $r=\rk(f_*\omega_{X/Y}^\nu)$ consider the family $f^{(r)}:X^{(r)} \to Y$, obtained as a desingularization of
the total space of the $r$ fold product 
$$
X^r=X\times_Y\cdots \times_Y X \>>> Y.
$$
Remark that 
$X^r$ is normal, Gorenstein with at most rational singularities. Hence flat base change implies that
$$
f^{(r)}_*\omega_{X^{(r)}/Y}^\nu=\bigotimes^rf_*\omega_{X/Y}^\nu.
$$
The family is perhaps not semistable, but replacing $Y$ by a covering sufficiently
ramified in $s\in S$ the pullback will have this property (This step is not needed if one replaces
in \ref{co.1} ``semistable'' by ``unipotent local monodromy''). So Proposition \ref{co.1} implies for
$\det(f_*\omega_{X/Y}^\nu)\subset f^{(r)}_*\omega_{X^{(r)}/Y}^\nu$ that
$$
\deg(f_*\omega_{X/Y}^\nu)\leq \frac{\rk(f_*\omega_{X/Y}^\nu)\cdot n\cdot \nu}{2}\cdot \deg(\Omega^1_Y(\log S)).
$$ 
\myqed 
\par\medskip \noindent
{\it Proof of \ref{co.1}.} \ 
For $\nu=1$ this is the inequality (\ref{eqhb.4}), obtained at the end of the last section.
Hence we will only consider the case $\nu >1$ in the sequel.

Assume that $\deg(\sH)> \frac{n\cdot \nu}{2}\cdot \deg(\Omega^1_Y(\log S))$. Replacing $Y$ by a finite covering, \'etale over
$U$, one may assume that 
$$
\deg(\sH) = \nu\cdot\rho > \frac{n\cdot \nu}{2}\cdot \deg(\Omega^1_Y(\log S)).
$$
Let $P$ be an effective divisor on $Y$ of degree $\rho$. Then $\sH\otimes \sO_Y(-\nu\cdot P)$
is in ${\rm Pic}_0(Y)$, hence divisible. So for some invertible sheaf $\sN$ of degree zero,
for $F=f^*P$ and for $\sL=\omega_{X/Y}\otimes f^*\sN\otimes \sO_X(-F)$ the sheaf
$\sL^\nu$ has a non-zero section $\sigma$. It gives rise to a cyclic covering of $X$ whose desingularization
will be denoted by $\hat{W}$ (see \cite{EV92}, for example). Then for some divisor
$\hat{T}$ the morphism $\hat{h}:\hat{W}\to Y$ will be smooth over $Y\setminus\hat{T}$, but not semistable. Choose $Y'$ to be a covering, sufficiently ramified, such that the 
pullback family has a semistable model over $Y'$.

Next choose $W'$ to be a $\Z/\nu$ equivariant desingularization of $\hat{W}\times_YY'$, and $Z$ to be a desingularization of the quotient. Finally let $W$ be the normalization of $Z$ in the function field of $\hat{W}\times_YY'$.
So we have a diagram
\begin{equation}\label{eqco.1}
\begin{CD}
W \>\tau >> Z \> \delta >> X' \> \varphi' >> X\\
\V h VV \V g VV \V f' VV \V f VV\\
Y' \> = >> Y' \> = >> Y' \> \varphi >> Y.
\end{CD}
\end{equation}
The $\nu$-th power of the sheaf $\sM=\delta^*\varphi'^*\sL$ has the section $\sigma'=\delta^*\varphi'^*(\sigma)$.
The sum of its zero locus and the singular fibres will become a normal crossing divisor after
a further blowing up. Replacing $Y'$ by a larger covering, one may assume that $Z\to Y'$ is semistable, 
and that $Z$ and $D$ satisfy the assumption iii) stated below. 

For a suitable choice of $T$ one has the following conditions:
\begin{enumerate}
\item[i.] $X'=X\times_YY'$, and $\tau:W\to Z$ is the finite covering obtained by taking the $\nu$-th root out of $\sigma' \in H^0(Z,\sM^\nu)$.
\item[ii.] $g$ and $h$ are both smooth over $Y'\setminus T$
for a divisor $T$ on $Y'$ containing $\varphi^{-1}(S)$. 
Moreover $g$ is semistable and the local monodromy of $R^nh_*\C_{W\setminus h^{-1}(T)}$ 
in $t\in T$ are unipotent. 
\item[iii.] $\delta$ is a modification, and $Z\to Y'$ is semistable. Writing $\Delta'=g^*T$ and $D$ for the zero divisor of $\sigma'$ on $Z$, the divisor
$\Delta'+D$ has normal crossing and $D_{\rm red}\to Y'$ is \'etale over $Y'\setminus T$.
\item[iv.] $\delta_*(\omega_{Z/Y'} \otimes \sM^{-1})=\varphi^*(\sN^{-1}\otimes \sO_Y(P))$,
\end{enumerate}
In fact, since $f:X\to Y$ is semistable, $X'$ has at most rational double points. Then 
$$
\delta_*(\omega_{Z/Y'} \otimes \delta^*\varphi'^*\omega_{X/Y}^{-1})=
\delta_*(\omega_{Z/Y'} \otimes \delta^*\omega_{X'/Y'}^{-1})=\delta_*\omega_{Z/X'}= \sO_{X'},
$$ 
which implies iv). The properties i), ii) and iii) hold by construction.

$W$ might be singular, but the sheaf $\Omega_{W/Y'}^p(\log \tau^*\Delta')=\tau^*\Omega_{Z/Y'}^1(\log
\Delta')$ is locally free and compatible with desingularizations.
The Galois group $\Z/\nu$ acts on the direct image sheaves
$\tau_*\Omega_{W/Y'}^p(\log \tau^*\Delta')$. As in \cite{EV92} or \cite[Section 3]{VZ05} one has the following description of the sheaf of eigenspaces.
\begin{claim}
Let $\Gamma'$ be the sum over all components of $D$, whose multiplicity
is not divisible by $\nu$. Then the sheaf
$$
\Omega^p_{Z/Y'}(\log (\Gamma'+\Delta'))\otimes \sM^{-1} \otimes \sO_{Z}\big(\big[\frac{D}{\nu}\big]\big),
$$
is a direct factor of ${\tau}_*\Omega^p_{W/Y'}(\log {\tau}^*\Delta')$. Moreover the $\Z/\nu$ action on $W$
induces a $\Z/\nu$ action on
$$
\W=R^nh_*\C_{W\setminus \tau^{-1}\Delta'}
$$
and on its Higgs bundle. One has a decomposition of $\W$ in a direct sum of sub variations of Hodge structures, given by the eigenspaces for this action, and the Higgs bundle of one of them is of the form
$ G=\bigoplus_{q=0}^n G^{n-q,q}$ for
$$
G^{p,q}=R^qg_*\big(\Omega^{p}_{Z/Y'}(\log (\Gamma'+\Delta'))\otimes \sM^{-1}
\otimes \sO_{Z}\big(\big[\frac{D}{\nu}\big]\big)\big).
$$
The Higgs field $\theta_{p,q}:G^{p,q} \to G^{p-1,q+1}\otimes \Omega^1_{Y'}(\log T)$ is induced by the edge
morphisms of the exact sequence
\begin{multline}\label{eqco.2}
0\>>>
\Omega^{p-1}_{Z/Y'}(\log (\Gamma'+\Delta'))\otimes {g}^* \Omega^1_{Y'}(\log T)\\
\>>> \Omega^{p}_{Z}(\log (\Gamma'+\Delta'))
\>>> \Omega^{p}_{Z/Y'}(\log (\Gamma'+\Delta')) \>>> 0,
\end{multline}
tensorized with $\sM^{-1} \otimes \sO_{Z}\big(\big[\frac{D}{\nu}\big]\big)$.\myqed  
\end{claim}
The sheaf 
$$
G^{n,0}=g_*\big(\Omega^n_{Z/Y'}(\log (\Gamma'+\Delta'))\otimes \sM^{-1}
\sO_{Z}\big(\big[\frac{D}{\nu}\big]\big)\big) 
$$
contains the invertible sheaf
$$
\sH=g_*\big(\Omega^n_{Z/Y'}(\log \Delta')\otimes \sM^{-1}\big)=
g_*(\omega_{Z/Y'} \otimes \sM^{-1})=\varphi^*(\sN\otimes\sO_Y(P)).
$$
of degree $\deg(\varphi)\cdot\rho$. 
\begin{claim}
$\sH$ generates a sub Higgs bundle
$$
(H=\bigoplus_{q=0}^nH^{n-q,q} ,\theta|_H)
$$ 
of $(G,\theta)$ where
$$
H^{n-q-1,q+1}={\rm Im}\big(\theta|_{H^{n-q,q}}:H^{n-q,q} \to G^{n-q+1,q+1}\otimes \omega_{Y'}(T)\big)\otimes \varphi^*\omega_{Y}(S)^{-1}.
$$
For some $q_0$ the sheaf $H^{n-q,q}$ is invertible of degree 
$$
\deg(\varphi)\cdot \big(\rho - q\cdot\deg(\Omega^1_Y(\log S))\big)
$$
for $q\leq q_0$ and zero for $q>q_0$.
\end{claim}
\begin{proof}
Writing $\Delta=f^*(S)$ consider the tautological exact sequences
\begin{equation}\label{eqco.3}
0\to
\Omega^{p-1}_{X/Y}(\log \Delta)\otimes {f}^* \Omega^1_{Y}(\log S)
\>>> \Omega^{p}_{X}(\log \Delta)
\>>> \Omega^{p}_{X/Y}(\log \Delta) \to 0,
\end{equation}
tensorized with 
$$
\omega_{X/Y}^{-1}=(\Omega^{n}_{X/Y}(\log \Delta))^{-1}.
$$
Taking the edge morphisms one obtains a Higgs bundle starting with the $(n,0)$ part $\sO_Y$. The sub 
Higgs bundle generated by $\sO_Y$ has $\Omega^1_Y(\log S)^{-q}$ in degree $(n-q,q)$.
Tensorizing with $\sN^{-1}\otimes\sO_Y(P)$ one obtains a Higgs bundle $H_0$ with
$$
H_0^{n-q,q}= \sN^{-1}\otimes \sO_Y(p)\otimes \omega_Y(S)^{-q}.
$$
On the other hand, the pullback of the exact sequence (\ref{eqco.3}) to $Z$ is a subsequence of
$$
0\to
\Omega^{p-1}_{Z/Y'}(\log \Delta')\otimes {g}^* \Omega^1_{Y'}(\log T)
\to \Omega^{p}_{Z}(\log \Delta')
\to \Omega^{p}_{Z/Y'}(\log \Delta') \to 0,
$$
hence of the sequence (\ref{eqco.2}), as well. So the Higgs field of $\varphi^*H_0$ 
is induced by the edge morphism of the exact sequence (\ref{eqco.2}), tensorized with
$$
\varphi'^*(\omega_{X/Y}^{-1}\otimes f^*(\sN^{-1}\otimes\sO_Y(P))),
$$ 
or with the larger sheaf $\sM^{-1} \otimes \sO_{Z}\big(\big[\frac{D}{\nu}\big]\big)$.

One obtains a morphism of Higgs bundles $\varphi^*H_0\to G$. By definition 
$$
\varphi^*H_0^{n,0} = \varphi^*(\sN^{-1}\otimes\sO_Y(P))= H^{n,0} \> \subset >> G^{n,0},
$$
and $H$ is the image of $\varphi^* H_0$ in $G$.
\myqed \end{proof}
The last claim implies that the degree of $H$ is  
\begin{multline*}
\deg(\varphi)\cdot \big((q_0+1)\cdot \rho - \sum_{q=0}^{q_0} q \cdot\deg (\Omega^1_Y(\log S))\big)=\\
\deg(\varphi)\cdot \big((q_0+1)\cdot \rho - \frac{q_0\cdot(q_0+1)}{2}\cdot\deg (\Omega^1_Y(\log S))\big)>\\
\deg(\varphi)\cdot \big(\frac{(q_0+1)\cdot n}{2} - \frac{q_0\cdot(q_0+1)}{2}\big)\cdot\deg (\Omega^1_Y(\log S)) \geq 0,
\end{multline*}
contradicting the polystability of $(G,\theta)$.
\myqed   

\begin{corollary}{\cite[Theorem 6]{VZ05}}\label{co.2}
Consider a semistable family $f:X\to Y$ of $n$-folds over a curve $Y$, smooth over $U=Y\setminus S$. Let $\V$ be a $\C$-sub variation of Hodge structures in
$R^nf_*\C_V$, with $(n,0)$ part $F^{n,0}$ and with a strictly maximal Higgs field.
Then the rational map $\pi:X\to \BP(F^{n,0})$,  defined by $f^*F^{n,0}\to \omega_{X/Y}$,
can not have a non-isotrivial image. In particular $\pi$ can not be birational.
\end{corollary}

\noindent{\it Sketch of the proof.} \ 
By Lemma \ref{in.7} one may assume that $F^{n,0}=\sL^n\otimes \sT$ for a logarithmic theta characteristic $\sL$ and for a unitary bundle $\sT$ on $Y$. 
Let $\sB_\nu$ denote the image of the multiplication map
$$
\mathfrak{m}:S^\nu(F^{n,0}) \>>> f_*\omega_{X/Y}^\nu.
$$
If the image of $\pi$ is non-isotrivial, the method used to prove \cite[Theorem 4.33]{V95} shows
that $\deg(\sB_\nu\otimes \sL^{-\nu\cdot n})$ has to be larger than zero, for some $\nu\gg 1$. 
For $r=\rk(\sB_\nu)$ one obtains an inclusion
$$
\det(\sB_\nu) \>>> \bigotimes^rf_*\omega_{X/Y}^\nu.
$$
Since 
$$
\deg(\sB_\nu) > r\cdot\nu\cdot n\cdot \deg(\sL)=\frac{r\cdot n\cdot\nu}{2}\cdot\deg(\omega_{X/Y}^\nu)
$$
this contradicts Proposition \ref{co.1}, applied to a desingularization of the $r$ fold
fibre product $X\times_Y\cdots \times_Y X \to Y$.
\myqed
\par\medskip \noindent
{\it Sketch of the proof of: ``\ref{in.13} $\Longrightarrow$ \ref{in.14}''.} \ 
We assume again that $\V$ has a maximal Higgs field, that $\rk(\V)>2$ and we write $F^{1,0}$
for the first Hodge bundle.

By Corollary \ref{co.2} it only remains to consider the case where the image $W$ of
$\pi$ is isotrivial. For $n=1$ one shows (replacing $Y$ by an \'etale covering) that $W\to Y$ is either of the form $C\times Y \to Y$, for some curve $C$ of positive genus, or it is a $\BP^1$ bundle. In the first case, 
it is easy to see, that $\BP(F^{1,0})$ is trivial, which implies by Lemma \ref{in.7} that $\V$ is the direct sum of copies of $\L$ contradicting Corollary \ref{in.13}.

A similar argument works in the second case, and we refer to \cite{VZ05} for the details.
\myqed   
\section{Maximal Higgs fields and Shimura curves}\label{mh}
In this section $f:X\to Y$ is a family of $g$-dimensional Abelian varieties. We assume as usual that the family is semistable, or slightly weaker, that
the local monodromies of $R^1f_*\C_V$ in $s\in S$ are all unipotent. We will sketch the proof
of Theorem \ref{in.8} under the quite restrictive assumption
\begin{assumption}\label{mh.1}
$\V=R^1f_*\C_V$ has no unitary part. 
\end{assumption}
Of course we will also assume that the Higgs field of $R^1f_*\C_V$ is strictly maximal.
So the Assumption \ref{mh.1} will exclude all non rigid families. The proof of Theorem \ref{in.8} without it is unfortunately much more difficult. 

For $S\neq\emptyset$ the Assumption \ref{mh.1} is not really a restriction.
If $R^1f_*\C_V=\V\oplus\U$ where $\V$ has a strictly maximal Higgs field, and where $\U$ is unitary, the assumption on $S$ implies by Lemma \ref{hb.8}, (2), that (replacing $U$ by an \'etale covering) this decomposition is defined over $\Q$ and that $\U$ is trivial. 
So $f:X\to Y$ is isogenous to the product of a family $f':X'\to Y$, with $\V$ as variation of Hodge structures, and a constant Abelian variety. In this case, the Assumption \ref{mh.1} holds if one replaces  $f:X\to Y$ by its ``moving part'' $f':X'\to Y$.

By Lemma \ref{in.7} there is a decomposition $R^1f_*\C_V=\L\otimes_\C\T$, and by the choice of
$\L$ one has $\det(\L)=\C$, hence $\det(\T)=\C$, as well. Again one can show, that such a decomposition
exists with $\L$ and $\T$ defined over real number fields. Here one applies the arguments
used in the proof of Lemma \ref{hb.8} to the local system $\E{\rm nd}(\L\otimes \T)$. 

One has a description of $g$-fold wedge products of tensor products (see \cite{FH91}, for example) in terms of partitions
$\lambda=\{\lambda_1, \ldots , \lambda_\nu\}$ of $g$. By definition
$\lambda_1, \ldots , \lambda_\nu$ are natural numbers with $g=\lambda_1 + \cdots + \lambda_\nu$. The partition $\lambda$ defines a Young diagram and a Schur functor ${\mathbb S}_\lambda$. Some standard elementary calculations show that:

\begin{lemma}\label{mh.2} \ 
\begin{enumerate}
\item[a.] If $k$ is odd, then for some partitions $\lambda_c$,
$$
\bigwedge^{k}(\L\otimes \T) = \bigoplus_{c=0}^{\frac{k-1}{2}}
S^{2c+1}(\L)\otimes {\mathbb S}_{\lambda_{2c}}(\T).
$$
\item[b.] If $k$ is even, then for some partitions $\lambda_c$,
$$
\bigwedge^{k}(\L\otimes \T) = S^{k}(\L)\oplus {\mathbb
S}_{\{2,\ldots,2\}}(\T) \oplus \bigoplus_{c=1}^{\frac{k}{2}-1}
S^{2c}(\L)\otimes {\mathbb S}_{\lambda_{2c}}(\T).
$$
\end{enumerate}
\end{lemma}

\begin{lemma}\label{mh.3} \
\begin{enumerate}
\item[a.] If $k$ is odd,
$$
H^0(Y,\bigwedge^k(\L\otimes\T))=0
$$
\item[b.] If $k$ is even, say $k=2c$, 
$$
H^0(Y,\bigwedge^k(\L\otimes\T))=H^0(Y,\bigwedge^k(\L\otimes\T))^{c,c}=
H^0(Y,{\mathbb S}_{\{2,\ldots,2\}}(\T)).
$$
\end{enumerate}
\end{lemma}
\begin{proof}
$S^\ell(\L)$ has a strictly maximal Higgs field for $\ell >0$, whereas for
all partitions $\lambda'$ the variation of Hodge structures
${\mathbb S}_{\lambda'}(\T)$ is again pure of bidegree $0,0$. 
The local system $S^\ell(\L)\otimes{\mathbb S}_{\lambda'}(\T)$ has again a strictly maximal Higgs field. A global section gives a trivial local subsystem, hence a Higgs subbundle
of the form $(\sO_Y,0)$, contradicting the strict maximality of the Higgs field.
So $S^\ell(\L)\otimes{\mathbb S}_{\lambda'}(\T)$ has no global sections. 

Hence $\bigwedge^k(\L\otimes \T)$ can only have global sections for $k$ even. In this case, the global sections lie in
$$
\det(\L)^c \otimes {\mathbb S}_{{\{2,\ldots,2\}}}(\T)={\mathbb S}_{{\{2,\ldots,2\}}}(\T).
$$
\myqed\end{proof}

Let $F$ be an Abelian variety and let $Q$ be the polarization, i.e. a non-degenerate
antisymmetric form on $H^1(F,\Q)$. The
Hodge group $\Hg(F)$ is defined in \cite{Mu66} (see also \cite{Mu69}) as the smallest $\Q$-algebraic subgroup of $\Sp(H^1(F,\Q),Q)$, whose extension to $\R$ contains the complex structure
$$
u:S^1 \>>> \Sp(H^1(F,\R),Q).
$$
Here $z$ acts on $(p,q)$ cycles by multiplication with $z^p\cdot \bar z^q$.

In a similar way, one defines the Mumford-Tate group $\MT(F)$.
The complex structure $u$ extends to a morphism of real algebraic groups
$$
h:{\rm Res}_{\C/\R}\G_m \>>>\Gl(H^1(F,\R)),
$$
and $\MT(F)$ is the smallest $\Q$-algebraic subgroup of $\Gl(H^1(F,\Q))$, whose extension to $\R$ contains
the image (see \cite{D72}, \cite{D82} \cite{Moo98} and \cite{Sch96}). Let us recall some of its properties,
stated in \cite{Moo98} and \cite{D72} with the necessary references. $\MT(F)$ is a reductive group,
and it preserves the intersection form $Q$ up to scalar multiplication.

$\MT(F)$ is also the largest $\Q$-algebraic subgroup of $\Gl(H^1(F,\Q))$,
which leaves all Hodge cycles of $F\times \cdots \times F$
invariant, hence all elements
$$
\eta\in H^{2p}(F\times \cdots \times F,\Q)^{p,p}
=\big[\bigwedge^{2p}(H^1(F,\Q)\oplus \cdots \oplus
H^1(F,\Q))\big]^{p,p}.
$$
For a smooth family of Abelian varieties $f:V \to U$ there exist the union $\Sigma$ of countably many
proper closed subvarieties of $Y$ such that $\MT(f^{-1}(y))$ is independent of $y$ for
$y\in U\setminus \Sigma$ (see \cite{D72}, \cite{Moo98} or \cite{Sch96}). Let us fix such a very general
point $y\in U\setminus \Sigma$ in the sequel and $F=f^{-1}(y)$. One defines $\MT(R^1f_*\Q_{V})$
to be the Mumford-Tate group of $F$.

Consider Hodge cycles $\eta$ on $F$ which remain Hodge cycles under
parallel transform. Then $\MT(R^1f_*\Q_{V})$ is the largest $\Q$-subgroup of $\Gl(H^1(F,\Q))$
which leaves all those Hodge cycles invariant (\cite[\S 7]{D82} or \cite[2.2]{Sch96}).
For an algebraic group $G$ the derived group $G^{\rm der}$ is the subgroup generated by all
commutators, or equivalently the kernel of the homomorphism $G\to G^{\rm ab}$ for $G$ to
the maximal Abelian quotient.

Let $\Mon^0$ be the algebraic monodromy group, i.e. the connected component of the Zariski closure of the image of the monodromy representation. Let us recall two results from \cite{D82} and \cite{A92} (see \cite[1.4]{Moo98}).

\begin{proposition}\label{mh.4} \
\begin{enumerate}
\item[a.] $\Mon^0$ is a normal subgroup of the derived subgroup $\MT(F)^{\rm der}$
of $\MT(F)$.
\item[b.] If for some $y'\in Y$ the fibre $f^{-1}(y')$ has complex multiplication,
then $\Mon^0=\MT(R^1f_*\Q_{V})^{\rm der}$.
\end{enumerate}
\end{proposition}
 Remark however, that up to now we do not know anything about the existence of points with complex multiplication for the families considered in Theorem \ref{in.8}. So instead of \ref{mh.4}, b), we will use: 
\begin{proposition}\label{mh.5}
Let $f:X\to Y$ be a family of $g$-dimensional Abelian varieties. Assume
that the local monodromies of $R^1f_*\C_V$ are all unipotent and that
its Higgs field is strictly maximal. Then
$\Mon^0=\MT(R^1f_*\Q_{V})^{\rm der}$.
\end{proposition}
\begin{proof}
Let us write $\MT=\MT(R^1f_*\Q_{V})$ and let us fix a very general fibre $F$ of $f$.
By \cite[4.4]{S92} $\Mon^0$ is reductive. By \cite[3.1 (c)]{D82} it is sufficient to show that each tensor
$$
\eta\in \bigwedge^{k}\big(H^1(F,\Q)\oplus \cdots \oplus
H^1(F,\Q)\big)= H^{k}(F\times \cdots \times F, \Q)
$$
which is invariant under $\Mon^0$ is also invariant under $\MT^{\rm der}$.
By abuse of notations, let us replace $F\times \cdots \times F$ by $F$.

A section $\eta\in H^{k}(F, \Q)\setminus \{0\}$ gives rise to a global section
$$
\eta\in H^0\big(U, \bigwedge^k \L\otimes_\C\T\big).
$$
By Lemma \ref{mh.3}, a) $k$ must be even, and by part b) the section $\eta$
is pure of bidegree $(\frac{k}{2},\frac{k}{2})$. So $\eta$ is a Hodge cycle,
and by definition it is invariant under $\MT$.
\myqed \end{proof}

Mumford defines in \cite{Mu66} a Shimura variety $\sX(\Hg,u)$ of Hodge type as a moduli scheme
of Abelian varieties (with a suitable level structure) whose Hodge group is
contained in $\Hg$. In \cite{Mu69} he gives an explicit construction. $\sX(\Hg,u)$ is the image of
$$
\Phi: \Hg_\R \>>> \Sp(H^1(F,\Q),Q)_\R/{\big( \begin{array}{c} \mbox{centralizer of the} \\ 
\mbox{complex structure } u \end{array}\big)} =\sH_g
$$
divided by an arithmetic subgroup $\Gamma \subset \Hg$.
The kernel of $\Phi$ is a maximal compact subgroup.

The monodromy group for $\sX(\Hg,u)$ is contained in $\Hg$, hence equal to
$$
\Hg^{\rm der}=\MT^{\rm der}\cap \Hg.
$$
Since
$$
(\MT^{\rm der}\times \MT) / \MT^{\rm der}
$$
is isogenous to $\MT$, we may replace in Mumford's construction $\Hg$ by $\Hg^{\rm der}$,
and the dimension of $\sX(\Hg,u)$ is the dimension of $\Phi(\Hg^{\rm der})$.
\par\medskip \noindent 
{\it Proof of Theorem \ref{in.8} under the Assumption \ref{mh.1}.}

The decomposition $R^1f_*\C_V=\L\otimes_\C\T$ implies that the representation
defining the local system has values in the product of $\Sl(2,\R)$
with a compact group. In particular, the dimension of $\sX(\Hg,u)$ is one,
and the induced morphism $\varphi':U\to \sX(\Hg,u)$ is dominant.

Consider the composite $\varphi$
$$
U\> \varphi'>> \sX(\Hg,u) \>>> \sA_g,
$$
where $\sA_g$ denotes a fine moduli scheme of Abelian varieties, with a suitable level structure. If $(F^{1,0}\oplus F^{0,1},\theta)$ denotes the
Higgs bundle of $R^1f_*\C_V$ the pullback of tangent sheaf of $\sA_g$ is given by
$\varphi^*T_{\sA_g}=S^2(F^{0,1})|_U$, and $T_U\to \varphi^*T_{\sA_g}$ is the map 
$$
T_U \>>> S^2(F^{0,1}|_U)\>>> (F^{0,1}\otimes {F^{1,0}}^\vee)|_U=\sL_0^{-1}\otimes\sT_0\otimes\sL_0^{-1}
\otimes \sT_0^\vee
$$ 
induced by $\theta$. Here $\sL_0$ is the restriction to $U$ of the logarithmic theta characteristic, and $\sT_0=\T\otimes
\sO_U$. So $T_U$ is just the subsheaf $\sL_0^{-2}\otimes \sO_U$, where
$\sO_U\subset \sH om (\sT_0,\sT_0)$ is given by the homotheties. In particular
$T_U$ is a direct factor of $\varphi^*T_{\sA_g}$, and $\varphi$ and $\varphi'$ are both \'etale.
\myqed
\par\medskip \noindent
{\it On the proof of Corollary \ref{in.10}.}
By Lemma \ref{hb.8}, (2), the assumption $S\neq \emptyset$ implies that the unitary system $\T$ becomes trivial over some \'etale covering, as well as the maximal unitary sub-local system $\U$ in $R^1f_*\C_V$. So we may assume, that this holds true on $Y$ itself. So $f:X\to Y$ is isogenous to a constant Abelian variety, corresponding to the unitary
subsystem $\U$ and the moving part $f':X'\to Y$. By abuse of notations, we will assume that $f=f'$, hence that $\U=0$. 

Assume first that the general fibre $F$ of $f$ is simple, hence $\U=0$, and $\T=\C^{\oplus g}$.

Sections in $H^0(U,\E{\rm nd}(R^1f_*\C_V))$ correspond to trivial Higgs subbundles
of $\sE nd ((\sL\oplus \sL^{-1})^{\oplus g})$. The rang of this bundle is $4\cdot g^2$ and a rank $3\cdot g^2$ subsystem has a strictly maximal Higgs field. The rank of the largest sub Higgs bundle with a trivial Higgs field is $g^2$, and it is concentrated in bidegree $(0,0)$. Since the corresponding local sub system of $\E{\rm nd}(R^1f_*\C_V)$ is unitary, Lemma \ref{hb.8} implies that it is defined over $\Q$
and that it trivializes over an \'etale covering.
So we may assume that
$$
H^0(U,\E{\rm nd}(R^1f_*\Q_V))=H^0(U,\E{\rm nd}(R^1f_*\Q_V))^{0,0}=\Q^{g^2},
$$
hence that ${\rm End}(F)_\Q=\Q^{g^2}$. Now take your favorite textbook on Abelian varieties
and look for simple ones, with ${\rm End}(F)_\Q=\Q^{g^2}$, and with a positive dimensional non-compact moduli space. 
The only examples you will find are families of elliptic curves.
So $f:X\to Y$ is a modular family of elliptic curves, completely determined by the local system $\L$.

If $F$ is not simple, apply the argument indicated above to each simple factor of the family.
\myqed
\par\medskip \noindent
{\it Sketch of the proof of Corollary \ref{in.11}.}\\
The description $R^1f_*\C_V=\L\otimes\T$ implies that
$$
H^0(U,\E{\rm nd}(R^1f_*\C_V))=H^0(U,\E{\rm nd}(R^1f_*\C_V))^{0,0}.
$$
Then the rigidity follows from \cite{F83}.
\myqed\par\medskip

Let $X_{\Hg}$ be the $\Hg(\R)^+$-conjugacy class in 
$$
{\rm Hom}_{{\rm alg. grp}/\R}(  {\rm Res}_{\C/\R}\G_m , \Hg_\R)
$$
containing $u \circ ({\rm Res}_{\C/\R}\G_m \to S^1)$. Here $+$
denotes the topological connected component.
For the reader's convenience we note that $(\Hg, X_\Hg)$ 
is a Shimura datum in the sense of \cite[2.1.1]{D79}: \\[.1cm]
$\Hg$ is reductive and \cite{D79}  Proposition 1.1.14 derives the
axioms (2.1.1.1) and (2.1.1.2) from the fact that 
$u$ defines a complex structure compatible with the
polarization. The axiom (2.1.1.3), i.e. the non-existence of a $\Q$-factor
in $\Hg^{{\rm ad}}$ onto which $h$ projects trivially,
follows from $\Hg$ being the smallest $\Q$-subgroup
of $\Sp(H^1(F,\Q),Q)$ containing $u$.
\par
\medskip
\noindent{\it Sketch of the proof of Theorem \ref{in.9}.}
Since $\Hg$ is reductive we may split the representation
$\Hg \to \Sp(H^1(F,\Q),Q)$ into a direct sum of irreducible
representations. We may split off unitary representations
over $\Q$ (\cite{K87} Proposition 4.11). 
\begin{claim}
For each of the
remaining irreducible representations there is an isogeny
$i: \Sl_2(\R) \times K \to \Hg_\R$, where $K$ is a compact
group, such that the composition
of $i$ with $\Hg_\R \to \Sp(H^1(F,\Q),Q)_\R$ is the tensor
product of a representation of $\Sl_2(\R)$ of weight one 
by a representation of $K$. 
\end{claim}
Assuming the Claim, consider a maximal compact
subgroup $K_1$ of $\Sl_2(\R) \times K$ that maps to the centralizer
of $u$ under  
$$
\Sl_2(\R) \times K \>>> \Sp(H^1(F,\Q),Q)_\R.
$$ 
The double quotient 
$$
\sX' = \Gamma' \setminus\!(\Sl_2(\R) \times K)/K_1
$$ 
is an unramified cover of $\sX(\Hg,U)$. Since a strictly
maximal Higgs field is characterized by the Arakelov
equality we may as well prove that the pullback
variation of Hodge structures has a strictly maximal Higgs field. Since the fundamental
group of $\sX'$ acts via $\Gamma'$ this follows immediately
from the claim and Lemma 2.1 in \cite{VZ04}.\par
\medskip\noindent
{\it Proof of the Claim.}
We first analyze $\Hg^{{\rm ad}}$ and the $\Hg^{{\rm der}}(\R)$-conjugacy
class $X_{\Hg}^{\rm ad}$ of maps ${\rm Res}_{\C/\R}\G_m \to
\Hg^{{\rm ad}}$ containing $(\Hg \to \Hg^{\rm ad}) \circ u$. 
Note that $X_{\Hg}^{\rm ad}$ is a connected component of $X_{\rm Hg}$. 
\newline
Each $\Q$-factor of $\Hg^{{\rm ad}}$ onto which $h$ projects non-trivially
contributes to the dimension of $\sH(\Hg,u)$. Since we
deal with Shimura curves $\Hg^{{\rm ad}}$ is $\Q$-simple
by (2.1.1.3). Let 
$$\Hg^{{\rm ad}}_\R = \prod_{i \in I} G_i$$
be its decomposition into simple factors. Then 
$X_{\Hg}^{\rm ad} = \prod X_i$ for $X_i$ a $G_i(\R)$-conjugacy
class of maps ${\rm Res}_{\C/\R}\G_m \to G_i$.
For the same reason, only one of the simple factors, say $G_1$, 
of $\Hg^{{\rm ad}}_\R = \prod_{i \in I} G_i$
is non-compact. The possible complexifications $(G_1)_\C$ are classified
by Dynkin diagrams. The property `Shimura curve', i.e.\ dimension one, implies 
that  $G_1 \cong {\rm PSl}(2,\R)$.
\par
Now we determine the possible representations. The universal
cover $\widetilde{G_1} \to G_1$ factors though 
$$
G:={\rm Ker}(\Hg \to \prod_{i \in I \setminus \{1\}} G_i)^0.
$$ 
We apply \cite{D79} Section 1.3 to 
$$ (G_1,X_{\Hg}) \leftarrow (G,X_{\Hg}) \to (\Sp,\sH). $$
Since a finite-dimensional representation of $G_1$ factors
through $\Sl(2,\R)$, we conclude that $G \cong Sl(2,\R)$.
Moreover, such a representation corresponds to a fundamental
weight, hence of weight one. Now we let $\widetilde{K}$ be
the universal cover of $\prod_{i\in I\setminus \{1\}} G_i$. 
Since $\Hg \to \Hg^{\rm ad}$ is an isogeny, there is a lift 
of the universal cover $\widetilde{K} \to \Hg$. This lift
factors though a quotient $K$ of $\widetilde{K}$ such that the 
natural map $\Sl(2,\R) \times K \to \Hg_\R$ is an isogeny. 
\par
Since we assumed the representation $\Hg_\R \to \Sp(H^1(F,\Q),Q)_\R$
to be irreducible, also 
$$
\rho: \Sl(2,\R) \times K \to \Hg_\R \to \Sp(H^1(F,\Q),Q)_\R
$$ 
is irreducible. Let $W \subset H^1(F,\R)$ be an irreducible 
(necessarily weight one) representation of $\Sl(2,\R) \times \{{\rm id}\}$.
Since $K$ is reductive, hence its representations are semisimple, 
$\rho$ is the tensor product of $W$ and the representation
${\rm Hom}_{\Sl(2,\R) \times \{{\rm id}\}}(W,H^1(F,\R))$ of $K$. 
This completes the proof of the claim.
\myqed

\section{Teichm\"uller curves and their variation of 
Hodge structures}\label{tc}

We start with the assumptions of Theorem \ref{in.12}: Let $f:X \to Y$ 
be a semistable family of curves of
genus $g \geq 2$, smooth over $U = Y \setminus S$ and let $V = f^{-1}(U)$.
Suppose that $R^1f_* \C_V$ contains a  sub variation of Hodge
structures $\L$ of rank two with a strictly maximal Higgs field. Let us for the
moment admit the following:
\begin{lemma} \label{LLoverRR}
$\L$ is defined over $\R$
\end{lemma} 
By Theorem \ref{in.2}  the universal covering of $U$ is the upper half plane $\sH$.
Let $\widetilde{\varphi}: \sH \to \sH$ the period map for $\L$.
The strictly maximality of the Higgs field is now equivalent to
$\widetilde{\varphi}$ being an isomorphism (compare \cite{VZ04} Lemma 2.1).
In different terms $\widetilde{\varphi}$ is an isometry for the Kobayashi metric
on $\sH$. The choice of a Teichm\"uller marking (see below) on one of the 
smooth fibres of $f$ defines a lift of the moduli map $m: U \to M_g$ to
Teichm\"uller space. We denote this lift by $\widetilde{m}: \sH \to T_g$
and we let $j: T_g \to \sH_g$ the natural map to the
Siegel half space, which associates to the curve its Jacobian and
to the Teichm\"uller marking  the choice of a symplectic basis.
Finally we let $p_{11}: \sH_g \to \sH$ the projection of the matrix
in $\sH_g$ to its $(1,1)$-entry. If we have chosen the symplectic basis $B$
suitably, i.e.\ such that the first pair in $B$ spans the fibres of $\L$,
then
$$ \widetilde{\varphi} = p_{11} \circ j \circ \widetilde{m}.$$
Since the composite map $\widetilde{\varphi}$ is a Kobayashi isometry, 
also  $\widetilde{m}$ is an isometry, if we provide the image with the 
restriction of the Kobayashi metric on $T_g$. By a theorem of Royden the 
Kobayashi metric on $T_g$ coincides with the Teichm\"uller metric. To sum up:
\newline
The image of $U$ in $M_g$ is an algebraic curve, whose lift in $T_g$ 
is a complex geodesic for the Teichm\"uller metric. These objects 
are called {\em Teichm\"uller curves}. We
have thus proved one implication of Theorem \ref{in.12}.1. 
\par
Before showing the converse implication we need some background on geodesics for
the Teichm\"uller metric. Details may be found in \cite{IT92}. 
\newline
One of the equivalent definitions of Teichm\"uller space is the following. 
Let $T_g$ be the space
of pairs $(R,h)$ of a Riemann surface plus a quasi-conformal mapping 
$h: R_0 \to R$ up to isotopy (the {\em Teichm\"uller marking}) 
from a `reference surface' $R_0$ to $R$.  
By Teichm\"uller's theorem the isotopy class of $h$ contains a unique
representative with minimal maximal dilatation. Its Beltrami
coefficient (a $(-1,1)$-form) $\mu = h_{\overline{z}}/h_z$ is of the form $\mu = t|q|/q$ 
for some quadratic differential
$q \in \Gamma(R_0,(\Omega^1_{R_0})^{\otimes 2})$ and $|t|<1$. Conversely, 
one can solve the Beltrami equation for any $\mu$ of with $|\mu|_\infty <1$. 
In particular for each $\mu = t|q|/q$ 
there is a quasi-conformal map $h: R_0 \to R$ for some Riemann surface $R$
whose Beltrami coefficient is $\mu$.
\par
A quadratic differential $q$ on a Riemann surface $R_0$ determines an atlas
of $R$ minus the zero set of $q_0$ such that the transition functions
are translations maybe composed by multiplication with $\pm 1$:
The charts $\psi_i: U_i \to \C$ are
given by integrating by a square root of $q$ that exists locally.
The condition for $h$ to have minimal maximal dilatation may
be rephrased by saying that $h$ is affine in the charts $\psi_i$ for
some $q$ on $R_0$ and some ('terminal') quadratic differential
on $R$. The Teichm\"uller geodesic from ${\rm id}$ to $h$ is given by the
conformal maps with Beltrami coefficient $ct|q|/q$ for $c \in [0,1]$.
Thus for each $q$, the set 
$$
\{\mu_q(t) := t|q|/q, \,|t| < 1\}
$$ 
defines totally geodesic subset of $T_g$.\par
We briefly recall the construction of the complex structure on $T_g$ that
makes $T_g \to M_g$ holomorphic. Given $(R,h)$ or rather its Beltrami
coefficient $\mu$ on $R_0$, pull it back to universal covering $\sH$ of $R_0$
and extend it by zero to obtain a Beltrami coefficient $\tilde{\mu}$ on $\C$. 
Solving the Beltrami equation for $\mu$ gives
a quasi-conformal mapping $\tilde{h}: \C \to \C$, conformal on 
the lower half plane $\sH^*$. The Schwarzian derivative of
$\tilde{h}|_{\sH^*}$ is a quadratic differential on $\sH^*$.
Moreover, since $\mu$ was pulled back from $R = \sH/\Gamma_R$, 
the Schwarzian derivative descends from $\sH^*$ to a quadratic differential 
$$
q^*_\mu \in H^0(\sH^*/\Gamma_R, (\Omega_{\sH^*/\Gamma_R}^1)^{\otimes 2}).
$$ 
The map $\mu \mapsto q^*_\mu$ factors through Teichm\"uller space and the 
complex structure on $H^0(\sH^*/\Gamma_R, 
(\Omega_{\sH^*/\Gamma_R}^1)^{\otimes 2})$ is the one we want.
\newline
The map $\mu \mapsto q^*_\mu$ has a section locally around $0$ (due to Bers)
given by $q^*(z) \mapsto -2 {\rm Im}(z)q^*(\bar{z})$. This implies that a
neighborhood of $0$ of 
$$
\{\mu_q(t), \,|t| < 1\}
$$ 
is a holomorphic
submanifold of $T_g$. Since the complex structure on $T_g$ is independent
of the base point, all of $\{\mu_q(t), \,|t| < 1\}$ is a holomorphic
submanifold.  
\par
To sum up: For a fixed Riemann surface $R_0$ and $q \in
\Gamma(R_0,(\Omega^1_{R_0})^{\otimes 2})$
fixed,  the set $t|q|/q$, $|t|<1$
is a totally geodesic and holomorphic submanifold of $T_g$, a {\em 
Teichm\"uller disc}.
\par
In order to characterize which (few!) {Teich}\-m\"uller discs descend to 
Teich\-m\"uller curves in $M_g$ 
we provide another description of a Teichm\"uller disc. There is a
natural $\Sl(2,\R)$-action on the bundle of triples 
$(R,\varphi,q)$ over $T_g$. Post-compose local charts $\psi_i$ of $R$ 
with $A \in \Sl(2,\R)$ acting on $\C \cong \R^2$. 
Since overlapping charts differ only by translations and multiplication
by $\pm 1$, the compositions
$$ 
\xymatrix{ \widetilde{\psi_i}: U_i \ar[r] & \C \cong \R^2 \ar[r]^A & \R^2 \cong \C  }
$$
determine a new complex structure on $R$. We let this new Riemann
surface be $A \cdot R$. There is a unique quadratic differential
$A\cdot q$ on $A \cdot R$ whose integration charts are $\widetilde{\psi_i}$.
\newline
The action of ${\rm SO(2,\R)}$ does not change the complex structure
of $R$. Hence the orbit $\Sl(2,\R)\cdot(R,\varphi,q)$ projects
to a disc in $T_g$. If we lift 
$$
\begin{array}{lcl}
\sH & \to & \Sl(2,\R) \\
\tau &\mapsto &\left(\begin{array}{cc} 1 & {\rm Re} \,\tau \\
0 & {\rm Im} \,\tau \end{array} \right) \\
\end{array}
$$
one verifies that the Beltrami coefficient of the composition 
$\widetilde{\psi_i}$ is $\frac{i-\tau}{i+\tau}\frac{|q|}{q}$. Hence the 
image of the $\Sl(2,\R)$-orbit is a Teichm\"uller disc.
\par
\begin{proposition} \label{affinegroup}
A Teichm\"uller disc descends to a Teichm\"uller curve, if
and only if the setwise stabilizer of the disc in the mapping
class group is a lattice in ${\rm Aut}(\sH) = \Sl(2,\R)$. This stabilizer 
coincides up to conjugation with a quotient (by a finite group)
of the group of orientation-preserving diffeomorphisms
that are affine in the charts provided by $q$.
\end{proposition}
\par
\noindent{\it Proof.} Proposition 3.2 and Proposition 3.3 in \cite{McM03}.
\myqed \medskip
\par
When referring to a Teichm\"uller curve in the sequel we always assume
that it is {\em generated by $q=\omega^2$ the square of a holomorphic one-form}.
If a Teichm\"uller curve is generated by a quadratic differential on $R$, 
take the double covering of $R$ where it admits a square root. The above
criterion gives an easy way to check that this pair generates again
a Teichm\"uller curve. 
\par \medskip
\noindent{\it Proof of Corollary \ref{in.13}.} This is \cite{McM03}
Theorem 4.2 in a different language. Suppose the contrary was true. 
Then $U \to M_g$ defines a Teichm\"uller
curve. Over the universal cover $\sH$ of $U$ we may choose sections
$\omega_j = \omega_j(\tau)$ for $j=1,2$ generating the $(1,0)$-part of $\L_j$.
We may moreover choose a symplectic basis $\{a_k,b_k\}$ for $k=1,\ldots g$ 
of $R^1 (f_\sH)_* \Z_V$ such that $\int_{b_k} \omega_j(\tau) = \delta_{j,k}$
for $\tau \in \sH$. We consider the maps 
$$
F_{j}: \sH \to \sH, \ \tau \mapsto \int_{a_j} \omega_j(\tau) \mbox{ \ \ for \ \ }j=1,2.
$$
Since the local systems 
$\L_j$ are strictly maximal Higgs both $F_1$ and $F_2$ are isometries. 
\newline
Let $R$ be the fibre of $f$ over $i \in \sH$.
Since $\sH \to T_g$ is the Teichm\"uller disc generated
by $\omega_1$ (say), Ahlfors' variational formula implies that
$$ \left.\frac{dF_j}{dt}\right\arrowvert_{t=i} = \int_R \omega_j(i) 
\frac{\overline{\omega_1(i)}}{\omega_1(i)}. $$
Hence the norm of $dF_j$ at $\tau = i$ in the hyperbolic metric  equals
$$ ||dF_j|| = \left\arrowvert\int_R \omega_j(i)
\frac{\overline{\omega_1(i)}}{\omega_1(i)}
\right\arrowvert\,/\,{\rm Im}
F_j(i) = \left\arrowvert\int_R \omega_j(i) \frac{\overline{\omega_1(i)}}{\omega_1(i)}\right\arrowvert\ 
\,/ \,\int_R|\omega_(i)|^2. $$ 
Since $F_j$ is an isometry, this implies by the Schwarz Lemma 
that  $||dF_j||= 1.$
Since $\omega_2$ is not proportional to $\omega_1$ this violates the
Cauchy--Schwarz inequality.
\myqed \medskip
\par
\noindent{\it Sketch of the proof of Lemma \ref{LLoverRR}.}
Since $\overline{\L}$ is also strictly maximal Higgs we have $\L \cong
\overline{\L}$ by Corollary \ref{in.13}. It remains to check that $\L$ can indeed be defined over $\R$,
see \cite{Moe2}, proof of Theorem 5.3, or repeat the argument used in the proof of Lemma \ref{hb.8}.
\myqed \medskip
\par
\noindent{\it Proof of Theorem\ \ref{in.12}.}
We fix a fibre $R$ of $f$ and $\omega \in \Gamma(R,\Omega^1_R)$ that
generates the Teichm\"uller curve.
By its construction as $\Sl(2,\R)$-orbit the fundamental group
of $U$ maps to a lattice $\Gamma$ in $\Sl(2,\R)$ and the subspace
$$
\langle {\rm Re}\, \omega, {\rm Im}\, \omega \rangle \subset H^1(R,\R)
$$ 
is invariant under the action of the fundamental group. Let $\L$ be the
corresponding irreducible rank two $\C$-local system. As in Lemma \ref{hb.8} one finds that $\L$ is defined over 
a number field, whose Galois closure we denote by $\tilde{K}$. Moreover
we let $K$ the number field generated by the traces of $\Gamma$. If we
apply Deligne's semisimplicity (Lemma \ref{hb.5}) 
to the polarized $\C$-variation of Hodge structures
$R^1 f_* \C_V$ the local system $\L$ and all its Galois conjugates
will appear. Moreover we may read the above argument of \cite{VZ04} Lemma
2.1 backwards to conclude that $\L$ is strictly maximal Higgs. This proves 1.
\par
By  Corollary \ref{in.13} the local system $R^1 f_* \C_V$ contains no other local subsystem 
isomorphic to $\L$. We claim that 
$$ R^1 f_* \C_V = \W \oplus \M, \quad \W_{\tilde{K}} = \bigoplus_{\sigma \in
\Gal(\tilde{K}/\Q)/\Gal(\tilde{K}/K)} \L^\sigma, $$
where the summands $\W$ and $\M$ are defined over $\Q$ and $\L^\sigma$
are pairwise non-isomorphic rank two local systems. The fact that
$L$ and $L^\sigma$ are isomorphic if and only if $\sigma$ fixes $K$
is not needed here, but see the remark at the end of this section. 
\par
We now prove 2. The $\Q$-sub variation $\W$ determines, after choosing a
$\Z$-lattice and fixing a polarization of type $\delta$, a family of Abelian
varieties $A \to U$ of dimension $r$. By \cite{F83} the tangent space to the
space of deformations of the moduli map $U \to A_{r,\delta}$ is a subspace of
the global sections of the local system ${\E}{\rm nd}(\W_\C)$ of bidegree
$(-1,1)$. Using above decomposition of $\W_\C$ into irreducible summands, one
checks (\cite{Moe2} Lemma 3.3) that this tangent space is trivial. 
Hence $U \to A_{r,\delta}$ is
rigid and $U$ is defined over a number field.
\myqed \vspace{.3cm}
\par
The above way of characterizing Teichm\"uller curves inverts the history of
these objects: The first examples, starting with Veech (\cite{Ve89}),
of Teichm\"uller curves were shown to have this property by
checking that the stabilizer of $\sH \hookrightarrow T_g$ is a lattice.
Next, in \cite{McM03}, Teichm\"uller curves were constructed as
the image in the moduli space $M_2$ of the intersection of two higher-dimensional 
$\Sl(2,\R)$-invariant
loci in the bundle of one-forms over $M_2$. Only
recently (\cite{BM05}) Theorem \ref{in.12} was used to construct
an infinite series of Teichm\"uller curves, including the examples
of Veech.
\par
\smallskip
We address once again the decomposition of the variation of Hodge
structures of a Teichm\"uller curves. Obviously, if $\L \cong \L^\sigma$
then $\sigma$ fixes the trace field. The converse is needed in
\ref{hasRM}. We reproduce the full proof, since it is, besides Corollary
\ref{in.13}, the only
argument in this context that relies essentially on Teichm\"uller theory.
\newline
By Proposition \ref{affinegroup}, elements in the fundamental group 
of $U$ may be represented by diffeomorphisms that are affine in the
$\omega$-charts. We want to show that if $\sigma$ does not fix the trace field, then
$\L \not\cong \L^\sigma$.
Choose a hyperbolic element $\gamma$ in the fundamental
group whose trace $t = {\rm tr}\gamma$ is not fixed by $\sigma$.
Lift $\gamma$ to a diffeomorphism $\phi$ of some smooth 
fibre $R$ of $f$. Let $\phi^*$ be the induced diffeomorphism
of $H^1(R,\R)$. The endomorphism $\phi^* + (\phi^*)^{-1}$
acts on $\L|_R$ by multiplication by $t = {\rm tr}\gamma$. 
We need to show that ${\rm Ker}(\phi^* + (\phi^*)^{-1} - t\cdot
{\rm id})$ has dimension two. It suffices to show that
the eigenvalues $\lambda^+$ (resp.\ $\lambda^-$) of $\phi^*$ of 
maximal (resp.\ minimal) absolute
values are unique and that these are the eigenvalues of $\phi^*$
acting on $\L|_R$. This is shown in \cite{McM03} Theorem 5.3: 
\par
The diffomorphism $\phi$ is not of finite order nor does it fix a closed
loop in $R$. Hence $\phi$ is pseudo-Ansov (\cite{FLP79}). For
such a diffeomorphism there are two transverse measured foliations 
$\mu^+$, $\mu^-$ on $R$, one is expanded by $\lambda^+$ and one is contracted
by $\lambda^- = (\lambda^+)^{-1}$. These two foliations represent
real cohomology classes in $\langle {\rm Re}\omega, {\rm Im}\omega \rangle
= \L|_R$. It remains to show that $\lambda^+$ is simple and of
largest absolute value. The corresponding statement for $\lambda^-$
follows by considering $\phi^{-1}$.
\par
The cohomology $H^1(R,\R)$ is spanned by $\R$-linear combinations
of the Poincar\'e duals $C^\vee$ to simple closed curves. We know
how $\varphi$ acts on $C^\vee$: Since $C$ is stretched in the
direction of $\mu^+$ and contracted in the direction of $\mu^-$, 
we have for $n \to \infty$ that $(\lambda^+)^{-n} (\phi^*)^n \,C^\vee
\to \alpha \mu^+$ for some $\alpha \in \R$. This proves the nonexistence
of an eigenspace for $\phi^*$ in $H^1(R,\R)$ with larger eigenvalue
than $\lambda^+$ and the simplicity of $\lambda^+$.
\par

\section{The only Teichm\"uller-Shimura curve}\label{to}

This section is a comparison between Teichm\"uller- and Shimura
curves. Recall that in the proof of Theorem \ref{in.12} we
defined the family Abelian varieties $A \to U$ corresponding to the local system $\W$. 
It is determined up to isogeny by the Teichm\"uller  curve and coincides with
the Jacobian of the universal family over the Teichm\"uller curve
in case $r=g$. 
\par
\begin{theorem} \label{hasRM}
The family $A \to U$ has real multiplication by $K$. The locus of
real multiplication is the smallest Shimura subvariety of $A_r$ that contains
the image of the moduli
map $U \to A_r$.
\end{theorem}
\par
\noindent{\it Proof.} For each $a \in K$ the cycle
$$ C_a := \bigoplus_{\Gal(\tilde{K}/\Q)/\Gal(\tilde{K}/K)} \sigma(a) \cdot 
{\rm id}_{\L^\sigma} \in {\rm End}(\W_{\tilde{K}}) = H^0(U,\E{\rm nd}(\W_{\tilde{K}})) $$
is defined over $\Q$ and of bidegree $(0,0)$. It is hence an
endomorphism of $A/U$ (\cite{D71} Remark 4.4.6). Since $K$ is a trace field
of a lattice in $\Sl(2,\R)$ it is real. By the classification 
of endomorphisms of Abelian varieties it is totally real. This
proves the first claim.
\par
For the second statement we have to show that the local systems
$\W^{\otimes m} \otimes (\W^\vee)^{\otimes m'}$ do not contain
Hodge cycles other than products and tensor powers of the cycles $C_a$. 
This is shown in \cite{Moe2} Lemma 3.3. The spirit of this Lemma
is similar to the decomposition into Schur functors of exterior
powers of the variation of Hodge structures over a Shimura curve in Section \ref{mh}.
The difference consists of an essential use of the uniqueness of 
the strictly maximal Higgs local system at some point.
\myqed \vspace{.3cm}
\par
We now address the question of classifying non-compact Shimura curves 
that lie entirely inside $M_g$ for $g \geq 2$. 
By Theorem \ref{in.8} the variation
of Hodge structures consists of unitary and strictly maximal Higgs
local subsystems. By Corollary \ref{in.14} the strictly maximal Higgs
part splits off a rank two local subsystem. By Theorem \ref{in.12}
such a curve is automatically also a Teichm\"uller curve. We may
now look at the decomposition of the variation of Hodge structures in the proof of 
Theorem \ref{in.12}. Since the Galois conjugate of a non-cocompact
lattice in $\Sl(2,\R)$ still contains non-trivial parabolic elements
none of the local systems $\L^\sigma$ is unitary. By Theorem \ref{in.8}
again there is no $\L^\sigma$ for $\sigma \neq {\rm id}$.
To sum up:
\par
\begin{lemma}
Suppose $U \to M_g$ is a non-compact Shimura curve. Then there 
is an unramified cover $U' \to U$ such that the pullback
of the universal family over $M_g$ to $U'$ has a family of
Jacobians ${\rm Jac(f)}:J \to U'$ with fixed part of dimension $g-1$, i.e.\ there
is an Abelian variety $A$ of dimension $g-1$ such that
$A \times U'$ injects into $J$.
\end{lemma} 
\par
\noindent{\it Proof.}
The only thing left to remark is that the unitary parts in the
variation of Hodge structures can be trivialized after an unramified cover $U' \to U$ (see Lemma \ref{hb.8}).
\myqed \vspace{.3cm}
\par
The possible dimensions $d$ of a fixed part in a family of Jacobians
of dimension $g$ over a one-dimensional base have been studied by Xiao. 
He proves in \cite{Xi87} that
$$  d \leq \frac{5g+1}{6}.$$
This implies that a fixed part of dimension $g-1$  can only occur
for $g \leq 7$.
Families of curves of genus $g=2$, $g=3$ and $g-4$ with 
$g-1$-dimensional fixed part are known to exist, but for
$g=5,6,7$ existence is still an open question. 
\par
Theorem \ref{in.15} may be considered as a non-existence statement
for such families under the additional assumption `Shimura curve'.
\par
\medskip
\noindent{\it Sketch of the proof of Theorem \ref{in.15}, b.}
There are four main ingredients: First, the large fixed part limits
the number of possible singular fibres the family $f: X \to Y$.
\newline
Second, the projection to $J \to J/A$ defines a covering $\pi: V \to E$. 
By \cite{Moe1} we may assume that $U=X(d)$ is a 
modular curve and $d = \deg(\pi)$. 
\newline
Third, all invariants, in particular some intersection numbers, 
of the complex surface $E$ (or rather its completion) are known.
We factor $\pi = i \circ \pi_{\rm opt}$, where $i$ is an isogeny
and $\pi_{\rm opt}$ does not factor through a non-trivial isogeny.
The knowledge of intersection number on $V$ will then be used to
bound the degree $d_{\rm opt}$ of $\pi_{\rm opt}$.
\newline
Finally, the translation structure on a fibre $R$ of $f$
given by the one-form $\omega \in H^0(R,\Omega^1_R)$ that generates
the Teichm\"uller curve is shown to have very special properties. 
In a case-by-case discussion depending on $d_{\rm opt}$ and the
number of zeros of $\omega$ we show that these properties are
absurd, except for the case of the Legendre family and
the family
$$f_{W}: y^4 = x(x-1)(x-t)$$
over $X(2)\cong \BP^1 \setminus \{0,1,\infty\}$ with local coordinate $t$
in genus three.
\par
We now give more details on the first step. 
Let $P_1,\ldots,P_s \in R$ be the zeros of
$\omega$ and let $k_i$ ($i=1,\ldots,s$) be the multiplicity of
the zero. Since $U \to M_g$ is a Teichm\"uller curve, we
may suppose by \cite{GuJu00} Theorem 5.5 that $\pi|_R$ is
ramified over one point $O \in E|\pi(R)$ only. Equivalently, 
the translation structure determined by $\omega$ is such that
$\pi|_{R \setminus \{P_1,\ldots,P_s \} }$ is a square-tiled covering
(see e.g.\ loc.~cit.). One can deduce from this that singular
fibres of $f$ do not contain separating nodes (\cite{Moe3} Proposition
2.3) and, more precisely, that the dual graph of a
singular fibre is a ring, say of length $\ell(y)$ for $y \in S = 
Y \setminus U$. (loc.\ cit.\ Lemma 2.1). 
\par
Concerning the second step, we recall that we constantly consider
Teich\-m\"uller and Shimura curves up to unramified cover outside $S$, 
since this does not change the property `strictly maximal Higgs'. 
We may hence replace $U$ by $X(d)$ for $d \geq 3$ sufficiently large
and post-compose $\pi$ by an isogeny to maintain the property
$d = \deg(\pi)$.
\par
Concerning the third step, we remark that the genus of $X(d)$,
the number of cusps and the self-intersection of the zero-section
of $E$ are well-known as functions of $d$. From the construction
of a Teichm\"uller curve as $\Sl(2,\R)$-orbit one deduces that
the zeros $P_i$ extend to sections $p_i: U \to V$ of $f$ for $i=1,\ldots,s$.
We denote the closure of the image of $p_i$ in $X$ by $\overline{P_i}$. Since
we assume $X$ to be a minimal semistable model, we have
$$ \omega_{X/Y} = \sO_X\left(\sum \overline{P_i}  + d \Delta_d \cdot R 
+ D\right),$$
where $R$ is a fibre of $f$, $D$ is a divisor supported in the
singular fibres of $f$ and 
$$\Delta_d = \frac{d^2}{24} \prod_{p|d} \left(1 - \frac{1}{d^2}\right). $$
Using the geometry of the singular fibres one checks (\cite{Moe3} Lemma 3.13)
that in fact $D=0$. From the intersection numbers on $E$ one deduces
that on $X$ we have 
$$\overline{P_i}= -d\Delta_d/(k_i+1).$$ Since $\L$
satisfies the Arakelov equality and since the remaining local sub-systems 
of $R^1 f_* \C$ are trivial, we have 
$$ \deg f_* \omega_{X/Y} = d \Delta_d.$$
\par
We claim that for each of the singular fibres $R_y$ for $y \in S$ of $f$ 
we have 
$$\Delta \chi_{\rm top} (R_y):= \chi_{\rm top} (R_y) - \chi_{\rm top} (R) = 
d/d_{\rm opt}.$$ 
Putting all these data in the Noether formula
$$ 12 \deg f_* \omega_{X/Y} - \sum_{y \in S} \Delta \chi_{\rm top} (R_y)
= \omega^2_{X/Y}$$
we obtain
$$ d_{\rm opt} = \frac{12}{\left(\sum_{i=1}^s \frac{k_1^2}{k_1 +1}\right)
+ 16 -4g}.$$
Since by definition $\sum_{i=1}^s k_i =2g-2$ there is only a finite
number of integer solution to this equation. We thus obtain the
desired bound on $d_{\rm opt}$.
\par
To prove the claim, notice first that $\Delta \chi_{\rm top} (R_y)$
coincides with the number of nodes of the dual graph of $R_y$ or
equivalently with the local component group of the N\'eron model
of the family of Jacobians around at $y$. The size of the component
group at $y$ is known to be $d$ for the universal family of elliptic
curves over $X(d)$. The details how to related the two groups 
are in \cite{Moe3} Proposition 3.15.
\par
The reader is referred to \cite{Moe3} for details on the last step.
\myqed \vspace{.3cm} 
\par
We remark that the exceptional family $f_W$ is generated by
the degree $8$ covering $\pi|_R$ of the once-punctured torus with Galois group
the quaternion group. The differential $\omega$ is the pullback
of a holomorphic differential on $E|_{\pi(R)}$. This covering $\pi|_R$ can
also be obtained by a double covering
of an elliptic curve ramified precisely over all four $2$-torsion points
post-composed by the multiplication by two.
\par



\end{document}